\documentclass[10pt,a4paper]{amsart}
\usepackage{a4}
\usepackage[utf8]{inputenc}
\usepackage[english]{babel}

\usepackage{amsfonts,amssymb,amsmath}
\usepackage[color=green!30]{todonotes}

\usepackage{xcolor}
\usepackage[colorlinks,
    linkcolor={red!50!black},
    citecolor={blue!50!black},
    urlcolor={blue!80!black}]{hyperref}
\newcommand{\mS}{\mathbb{S}}
\newcommand{\mC}{\mathbb{C}}

\newcommand{\mP}{\mathbb{P}}
\newcommand{\mT}{\mathbb{T}}
\newcommand{\mZ}{\mathbb{Z}}
\newcommand{\mK}{\mathbb{K}}
\newcommand{\cO}{\mathcal{O}}

\newcommand{\Si}{\Sigma}

\newcommand{\bm}{\begin{pmatrix}}
\newcommand{\ema}{\end{pmatrix}}
\newcommand{\bsm}{\left(\begin{smallmatrix}}
\newcommand{\esm}{\end{smallmatrix}\right)}

\newcommand{\al}{\alpha}
\newcommand{\be}{\beta}
\newcommand{\g}{\gamma}
\newcommand{\G}{\Gamma}
\newcommand{\de}{\delta}

\newcommand{\la}{\lambda}

\newcommand{\benu}{\begin{enumerate}}
\newcommand{\eenu}{\end{enumerate}}
\newcommand{\brho}{{\bar{\rho}}}
\newcommand{\Adr}{\Ad\circ\rho}
\newcommand{\Adbr}{\Ad\circ\bar{\rho}}
\newcommand{\Ada}{\Ad\circ\alpha}

\usepackage{tikz}
\usetikzlibrary{matrix}
\usepackage[all]{xy}
\usetikzlibrary{calc, matrix, arrows, cd}
\tikzset{node distance=2cm, auto}

\newcommand{\Ad}{\operatorname{Ad}}
\newcommand{\Aut}{\operatorname{Aut}}
\newcommand{\SL}{\operatorname{SL}}
\newcommand{\SU}{\operatorname{SU}}
\newcommand{\GL}{\operatorname{GL}}
\newcommand{\Hom}{\operatorname{Hom}}
\newcommand{\Ext}{\operatorname{Ext}}

\newcommand{\tor}{\operatorname{tor}}
\newcommand{\Spec}{\operatorname{Spec}}
\newcommand{\Tr}{\operatorname{Tr}}
\newcommand{\ev}{\operatorname{ev}}
\newcommand{\im}{\operatorname{im}}
\newcommand{\coker}{\operatorname{coker}}
\newcommand{\rk}{\operatorname{rk}}
\newcommand{\slf}{\mathfrak{sl}}
\newcommand{\mfp}{\mathfrak{p}}

\newtheorem{theorem}{Theorem}[section]

\newtheorem{corollary}[theorem]{Corollary}
\newtheorem{lemma}[theorem]{Lemma}
\newtheorem{proposition}[theorem]{Proposition}

\newtheorem*{question}{Question}

\newtheorem*{claim}{Claim}
\newtheorem*{notation}{Notation}

\theoremstyle{definition}
\newtheorem{definition}[theorem]{Definition}

\newtheorem{example}[theorem]{Example}

\newtheorem{remark}[theorem]{Remark}

\begin{document}

\title{Reidemeister torsion form on character varieties}
\begin{abstract}
In this paper we define the adjoint Reidemeister torsion as a differential form on the character variety of a compact oriented 3-manifold with toral boundary, and prove it defines a rational volume form. Then we show that the torsion form has poles only at singular points of the character variety. In fact, if the singular point corresponds to a reducible character, we show that the torsion has no pole under a generic hypothesis on the Alexander polynomial, else we relate the order of the pole with the type of singularity. Finally we consider the ideal points added after compactification of the character variety. We bound the vanishing order of the torsion by the Euler characteristic of an essential surface associated to the ideal point by the Culler-Shalen theory. As a corollary we obtain an unexpected relation between the topology of those surfaces and the topology of the character variety.
\end{abstract}
\author{Leo Benard}
\address{Mathematisches Institut, Georg-August Unversit\"at, G\"ottingen, Deutschland}
\email{leo.benard@mathematik.uni-goettingen.de}
\date{}

\maketitle

\section{Introduction}
The Reidemeister torsion is a combinatorial invariant of a (co-)homological complex, celebrated in the 30's for being able to distinguish non-homeomorphic lens spaces, and for completing their classification. One of the main feature of the torsion is that it turns out to be a topological invariant of manifolds. This deep result is know as Chapman--Cohen's theorem, see \cite{Chapman, Cohen}. Specifically, the Reidemeister torsion is a topological invariant $\tor(M,\rho)$ where $M$ is a 3-dimensional manifold and $\rho$ is a representation of its fundamental group $\pi_1(M)$ into a Lie group $G$. If the twisted cohomology groups $H^i(M,\rho)$ vanish for all $i \ge0$, then the torsion is a numerical invariant defined up to sign. In the general case, we may interpret it as a volume element in the twisted cohomology, that is by definition an element in the one-dimensional vector space
\begin{equation*}
\det(H^*(M,\rho))=\bigotimes_{i=0}^3\det(H^i(M, \rho)) ^{(-1)^{i}}.
\end{equation*}

Moreover, if $\rho$ and $\rho'$ are conjugated representations, then there is a canonical isomorphism $\det(H^*(M,\rho))\simeq \det(H^*(M,{\rho'}))$ that preserves the torsion. Hence it is natural to define the Reidemeister torsion as a section of some line bundle over the character variety.

We will be interested in the case where $M$ is a 3-manifold with toral boundary (e.g. a knot exterior). Its fundamental group acts on the Lie algebra $\slf_2(\mC)$ of the group $\SL_2(\mC)$ by composition of the adjoint action with any representation $\rho\colon\pi_1(M) \to \SL_2(\mC)$, and gives rise to the \emph{adjoint torsion}. Joan Porti in his Phd thesis \cite{Porti97} defined the adjoint torsion as an analytic function on a Zariski open subset of the character variety depending on a choice of a boundary curve. Many computations have been performed by J. Dubois, Vu Huynh, Yoshikazu Yamaguchi in \cite{Dub06,DHY09} and the torsion has been extended to the whole character variety by Dubois--Garoufalidis in \cite{DG09}. The torsion of the complex induced by the standard action on $\mC^2$ has been studied among others by Kitano in \cite{kitano1994,Kit94}, see also \cite{Ben17}.
In this article we will consider the adjoint torsion as a rational volume form on the character variety. This point of view has a long story, initiated by Johnson in his unpublished notes. Then Witten (\cite{Witten}), Dubois (\cite{Duboisthesis}), Park (\cite{Park}) and Frohman--Kania-Bartoszynska (\cite{Frohman-KB}) have developed the theory of the Reidemeister torsion as a volume form, essentially on the $\SU(2)$-character variety. In this article we will follow the approach initiated by Julien March\'e in~\cite{Mar15}.
More precisely, if the boundary of $M$ is a torus, the torsion is a rational volume form on the {\it augmented character variety} which is the following $2$-fold covering of the character variety:
\begin{equation*}
\bar{X}(M)=\{(\rho\colon\pi_1(M)\to\SL_2(\mC),\lambda\colon\pi_1(\partial M)\to \mC^*), \Tr \rho|_{\pi_1(\partial M)}=\lambda+\lambda^{-1}\}/\!/\SL_2(\mC).
\end{equation*}

In this paper we assume that $X(M)$ is 1-dimensional and \textit{scheme reduced} (in fact, slightly weaker hypotheses will be sufficient to our purpose). The first assumption is guaranteed by the assumption that $M$ is \textit{small}, that is without closed essential oriented surfaces not parallel to the boundary.

Let $\bar X$ be an irreducible component of $\bar{X}(M)$ containing the character of an irreducible representation and let $Y$ be its smooth projective model. It is a smooth compact curve obtained from $\bar X$ by desingularizing and adding a finite number of points at infinity: we call those points \textit{ideal points} of $Y$ and the others are called \textit{finite points}. We will denote by $v$ an element of $Y$ that we may view as a valuation on the function field $\mC(Y)\simeq\mC(\bar X)$, with local ring at $v$ denoted by $\cO_v$. The adjoint torsion will be denoted by $\tor(M, \Adr)$ and seen as an element of $\Omega_{\mC(Y)/\mC}$. 
The first result in this article is the following theorem:

\begin{theorem}
\label{theo:IntroFinite1}
Let $v$ be a finite point of $Y$, then $tor(M, \Adr)$ does not vanish $v$. In addition, 
if $v$ projects to an irreducible character in $X(M)$, then the order of the pole of $\tor(M)$ at $v$ is the length of the torsion part of the module $\Omega_{\mC[\bar{X}(M)]/\mC}\otimes \cO_v$. 
In particular if $v$ projects to a smooth point of $X(M)$ then the torsion has no pole at $v$. 
\end{theorem}
The length of the torsion part of the module $\Omega_{\mC[\bar{X}(M)]/\mC}\otimes \cO_v$ is an invariant of branch of the singularity which can be computed explicitly, as we will explain in Section \ref{sub:Singularities}.

Assume that $M$ is a rational homology circle, and let $\varphi\colon\pi_1(M) \to \mZ$ be the abelianization homomorphism. If $v$ projects to a reducible character $\chi$ in the character variety, this character can be written $\chi = \la^\varphi+\la^{-\varphi}$, for some $\la \in \mC^*$. Then it is a well-known result (due to Burde \cite{Bur} and to de Rham \cite{DR67} independently) that in this case $\la^2$ turns  to be a root of the Alexander polynomial $\Delta_M$ of $M$. We prove:
\begin{theorem}
\label{theo:IntroFinite2}
If $v$ projects to a reducible character in $X(M)$ such that $\la^2$ is a root of order one of $\Delta_M$, then the torsion has no pole at $v$.
\end{theorem}
If $v$ is an ideal point of $Y$, then the Culler-Shalen theory associates to $v$ an action of $\pi_1(M)$ on the Bass-Serre tree of $\SL_2(\cO_v)$ which itself yields an essential surface, denoted by $\Sigma$, in $M$. We say that such a surface $\Sigma$ is dual to the ideal point $v$. 

\begin{theorem}
\label{theo:IntroIdeal}
Let $v$ be an ideal point of the augmented character variety $Y$ and $\Sigma$ be an essential separating surface dual to $v$. We suppose that $\Sigma$ is a union of parallel connected copies $\Sigma_1 \cup \ldots \cup \Si_n$ and that both components of $M \setminus \Sigma_i$ are handlebodies. Let us also assume that $Y$ contains the character of a representation whose restriction to $\Sigma$ is irreducible. 
Then the torsion $\tor(M, \Adr)$ has vanishing order at $v$ bounded above by $-n(\chi(\Sigma)+1)$.
\end{theorem}
We say that a surface $S$ in $M$ is {\it free} if its complement is a union of handlebodies. Many natural constructions yield such surfaces. For example, take a knot diagram and consider the checkerboard surfaces (for an example of such a surface, see the left-hand side of Figure \ref{noeud}). If one of them, say $\Sigma$, is an essential non orientable surface in $M$, then the boundary of a tubular neighborhood of $\Sigma$ is orientable, remains incrompressible and does split $M$ into two handlebodies, as can be easily seen (both part of its complement retract onto a graph). In fact, any essential surfaces is free when $M$ is small. On the other hand, any essential surface whose class in $H_2(M, \partial M)$ is non zero will be non-separating in $M$. 

\medskip

We deduce from this theorem an unexpected relation between the genus of the character variety of $M$ and the genus of the essential surfaces in $M$. 
More precisely, suppose that $M$ is a knot complement whose character variety is one dimensional. Pick a smooth component of the variety, and assume that each ideal point $y$ in its smooth projective model  $Y$ corresponds to an essential surface $\Sigma_y$ that satisfies the hypothesis of Theorem \ref{theo:IntroIdeal}. Let us further assume that the Alexander polynomial of $M$ has only simple roots. Then 
\begin{equation*}
\label{equa:Intro}
\chi(Y) \ge \sum\limits_y n_y(\chi(\Sigma_y)+1)
\end{equation*}
In the simple case where the surfaces $\Sigma$ are connected, it turns into 
\begin{equation}
\label{equa:Topology}
\chi(Y) \ge \sum_y (\chi(\Sigma_y) +1).
\end{equation}

\begin{example}
It is shown in Hatcher--Thurston  \cite{HT85} that the knot $5_2$ has two separating essential surfaces in its complement: $\Sigma_1$ whose Euler characteristic is $-4$, and $\Sigma_2$ whose Euler characteristic is $-2$ (see figure \ref{noeud}). The third essential surface is the minimal Seifert surface $S$.

\begin{figure}[h]
\begin{center}
\def\svgwidth{0.7\columnwidth}
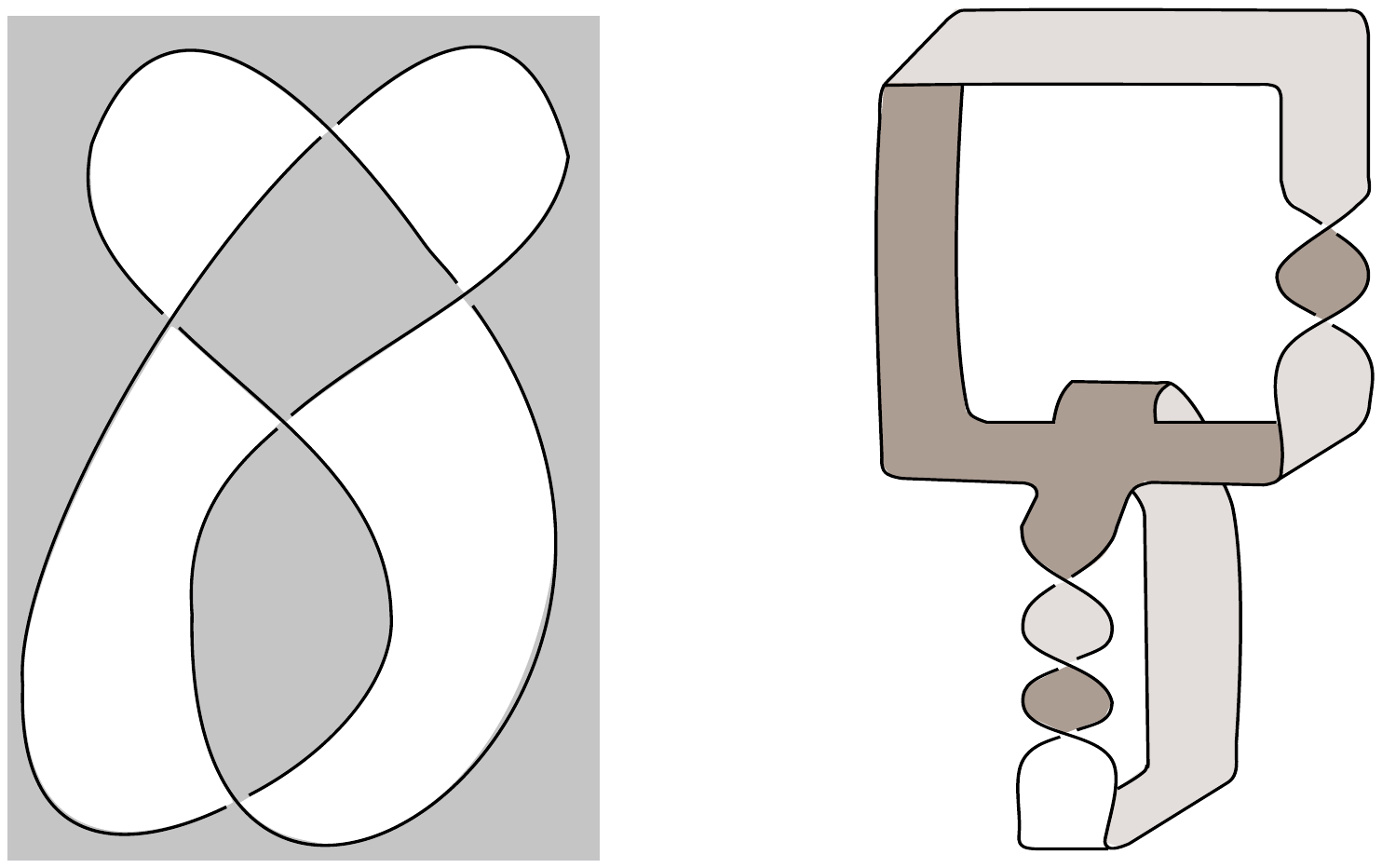
\caption{\label{noeud}Incompressible surfaces in the complement of the knot $5_2$. The surface $\Sigma_1$ is the orientation covering of the non-oriented surface colored on the left. It can be thought as the boundary of a tubular neighborhood of this non orientable surface. The surface $\Sigma_2$ is the orientation covering of the surface colored on the right. It can be obtained as follows: consider two parallel copies of each twisted bands above and below the square in the middle, and plumb them along this square. The result is connected because the bands below have an odd number of twists, and this is our surface $\Sigma_2$}
\end{center}
\end{figure}

The (geometric component of the) character variety $X$ has 3 ideal points, each of them corresponds to one of the essential surfaces described above. The torsion vanishes at order $1$ on the ideal point corresponding to  $\Sigma_2$, and at order $3$ at the ideal point corresponding to $\Sigma_1$. The Pl\"ucker formula yields $\chi(X)=-2$, and the covering map $Y \to X(M)$ ramifies at six points, hence $-\chi(Y)=10$ by the Riemann-Hurwitz formula. Since the torsion has zeros only at the three ideal points, and no poles, a direct computation using (\ref{equa:Topology}) shows that it must vanish at order $1$ at the ideal point corresponding to the Seifert surface.
\end{example}

\begin{question} 
Is the bound of Theorem \ref{theo:IntroIdeal} always sharp?
\end{question}
In all the examples we have listed in Section \ref{sec:Examples} for connected essential surfaces, it happens to be an equality. A careful examination of the proof shows that it has to be generically the case. The lack of equality should be interpreted as a non-transversal situation.

%

\medbreak 

Throughout this paper $k$ will denote an algebraically closed field of characteristic 0. The paper is organized as follows:
in Section \ref{sec:Character} we introduce character varieties and various related notions, in Section \ref{sec:Twisted} we define the vector space of rational differential forms on the character variety and relate it with the twisted cohomology of $M$, and we give a quick survey of the Culler-Shalen theory, and in Section \ref{sec:Reidemeister} we define the Reidemeister torsion form. In Section \ref{sec:Finite} we prove Theorems \ref{theo:IntroFinite1} and \ref{theo:IntroFinite2}, in Section \ref{sec:Examples} we relate the torsion form with previous works, and compute it explicitly on several examples, and finally in Section \ref{sec:Ideal} we prove Theorem \ref{theo:IntroIdeal}.

\subsection*{Acknowledgement}
This work has been conducted during the PhD thesis of the author, hosted by the Institut de Math\'ematiques de Jussieu in Sorbonne Universit\'e. The author thanks his advisor Julien March\'e for his inestimable time and help. He also thanks  Michael Heusener and Joan Porti for many helpful discussions that have widely contributed to improve this article, and Jinsung Park and Seokbeom Yoon for discussions that led to find a mistake in a previous version of Theorem 0.1, where the role of poles and zeros were inverted. Finally, the author is indebted to the anonymous referees for their very profitable remarks and suggestions.

\section{Character varieties and tautological representation}
\label{sec:Character}
In this section we introduce the basic material of this article: in Subsections \ref{sub:Rep} and \ref{sub:Char} we define the character variety of a finitely generated group, and in Subsection \ref{sub:Characters} we discuss various properties of characters. In this article we will deal with a refinement of the character variety that we define in Subsection \ref{sub:AugmentedVariety}. Then we define a crucial tool of the construction of the torsion form: the tautological representation in Subsection \ref{sub:Tautological}. We end the section with several examples in Subsection \ref{sub:Examples}

\subsection{Representation varieties}
\label{sub:Rep}
In this subsection we define the representation variety of a finitely generated group into  $\SL_2(k)$.

\medskip

Let $k$ be an algebraically closed field of characteristic 0, and let $\G$ be a finitely generated group. In the sequel the group $\G$ will always be the fundamental group of a 2- or 3-manifold, but we give definitions in this more general setting. For this subsection we fix $S = \lbrace \g_1, \ldots, \g_n \rbrace$ a generating set for $\G$.

We say that a \textit{representation} is a group homomorphism $\al \colon \G \to \SL_2(k)$. The \textit{representation variety} is the set $R(\G) = \Hom(\G,\SL_2(k)) =\lbrace\al\colon\G\to\SL_2(k) \rbrace$. The map
$$\iota_S \colon R(\G) \hookrightarrow \SL_2(k)^n$$
$$\al \mapsto (\al(\g_1), \ldots, \al(\g_n))$$
endows $R(\G)$ with a structure of an algebraic set whose image is defined as the zero locus in $\SL_2(k)^n$ of a finite set of polynomials given by the group relations. 

The algebra of functions of the representation variety is $$k[R(\G)]= k[X^{i,j}_\g, 1 \le i,j \le 2, \g \in \G]/(X_e -I, X_\g X_\de -X_{\g\de},\g, \de \in \G)$$ where for any $\g \in \G, X_\g$ denotes the matrix $\bsm X^{1,1}_\g&X^{1,2}_\g\\X^{2,1}_\g & X^{2,2}_\g \esm$ ($I$ denotes the matrix $\bsm 1&0\\0&1\esm$). This algebra is finitely generated, since any choice of a generating set provides a finite set of generators.

\subsection{Character varieties}
\label{sub:Char}
In this subsection we define the $\SL_2(k)$-character variety of a finitely generated group. It is classically defined as an algebro-geometric quotient, and we give an equivalent definition at the end of the subsection.

\medskip

The algebraic group $\SL_2(k)$ acts by conjugation on the representation variety~$R(\G)$. Two representations $\al, \al' \colon \G \to \SL_2(k)$ are conjugate if there exists $M \in \SL_2(k)$ such that for every $\g$~in $\G$, the matrix $\al(\g)$ is equal to $M\al'(\g)M^{-1}$. 

The action of the group $\SL_2(k)$ on $R(\G)$ induces a natural action on its algebra of functions $k[R(\G)]$ by pre-composition. The sub-algebra of invariant functions is 
\begin{equation*}
k[R(\G)]^{\SL_2} = \lbrace P\in k[R(\G)] \vert \ M\cdot P = P \text{ for all } M \in \SL_2(k) \rbrace.
\end{equation*} 
It is known to be finitely generated, although this is a delicate problem first answered by Hilbert (\cite{Hilbert}) in the late nineteenth century. There is an amount of good references on the topic, let us just mention \cite{MFK} and \cite{KP96}.

Recall that any finitely generated $k$-algebra $A$ is the quotient of a polynomial algebra, in other words there is an exact sequence 
$0 \to I \to k[X_1, \ldots, X_n] \to A \to 0$. By Hilbert's basis theorem, the ideal $I$ is finitely generated, hence $A$ defines an algebraic set $V(I) \subset k^n$, namely the zero-locus of any generating set of polynomials for $I$. Up to isomorphism, this set does not depend on the presentation of $A$, hence we denote by $\Spec(A)$ the affine algebraic variety defined by $A$.
\begin{remark}
\label{rem:reduced}
In general, the use of the term variety is reserved to irreducible and reduced algebraic sets.
An \textit{irreducible} set is a set which is not a reunion of two proper closed subsets. An \textit{irreducible component} is a maximal irreducible subset. Given a ring $R$, its spectrum $\Spec(R)$ is said to be \textit{reduced} if $R$ does not contain any nilpotent element. In particular an irreducible component is reduced. We will call many algebraic sets varieties despite they have no reason to be irreducible, nor reduced.
\end{remark}
\begin{definition}
The \textit{character variety} $X(\G) = R(\G) /\! / \SL_2(k)$ is the spectrum $\Spec(k[R(\G)^{\SL_2(k)}])$ of the sub-algebra of invariant functions.
\end{definition}
It is usually called the algebro-geometric quotient of $R(\G)$ by $\SL_2(k)$. Let us list without proof some of its properties:
\begin{itemize}
\item It comes with a projection map $\pi\colon R(\G) \to X(\G)$ that satisfies the following universal property: for any $\SL_2(k)$-invariant morphism $F\colon R(\G) \to Y$, with $Y$ an algebraic variety, there is a unique map $F'\colon X(\G) \to Y$ such that $F = F' \circ \pi$.
\item The $k$-points of this quotient are in bijection with the closed orbits of  $\SL_2(k)$ acting on $R(\G)$, or with conjugacy classes of semi-simple (or completely reducible) representations of $\G$ into $\SL_2(k)$. In other words, orbits whose closure intersect in $R(\G)$ are identified in $X(\G)$.
\item It is the biggest Hausdorff quotient of the topological quotient $R(\G)/\SL_2(k)$.
\end{itemize}

The following functions will play the role of coordinate functions on the character variety. Those functions are involved in the classical theory of character varieties of 3-manifolds, see \cite{CS83, CS84, Sha02}.
\begin{definition}
For any $\g \in \G$, we define the \textit{trace function} $I_\g: R(\G) \to k$ by $I_\g(\al) =~\Tr(\al(\g))$. Those functions are invariant under the action of $\SL_2(k)$. By the first property above about algebro-geometric quotients, they define functions on the quotient $X(\G)$, that we still denote by~$I_\g$.
\end{definition}
The following lemma is straightforward, but crucial:
\begin{lemma}
\label{lem:TraceRel}
For any $\g, \de \in \G$, the identity $I_\g I_\de = I_{\g\de}+I_{\g\de^{-1}}$ holds on $R(\G)$. 
\end{lemma}

This lemma motivates the definition of the trace algebra 
\begin{equation*}
B[\G] = k[I_\g, \g \in \G]/(I_e-2,I_\g I_\de - I_{\g\de}-I_{\g\de^{-1}}, \g, \de \in \G).
\end{equation*}
It turns out that this algebra is finitely generated, there is even an explicit bound on the number of generators in \cite[Proposition 1.4.1]{CS83}. A deep theorem of invariant theory states that this algebra is isomorphic to the algebra of invariants $k[X(\Gamma)]$ (\cite{P87,PS00}, see also \cite[Theorem 2.15]{Mar15}).

\subsection{Irreducible and reducible characters}
\label{sub:Characters}
In this subsection we define reducible and irreducible characters. Then we state a theorem due to Kyoji Saito, that will be used in particular in Subsection~\ref{sub:Tautological} to construct the tautological representation.

\medskip

It appears in Subsection \ref{sub:Char} that the character variety is defined as a scheme: it is the spectrum of the algebra of invariants $k[R(\G)]^{\SL_2(k)}$. Hence we can consider points of the character variety in the sense of algebraic geometry: given an integral $k$-algebra $R$, an $R$-point of $X(\G)$ is a class of $k$-algebras morphism $\chi \colon k[X(\G)] \to R$, where $\chi$ is identified with $\chi'$ if they have the same kernel. In particular a $k$-point is a maximal ideal $\mathfrak{m}$ given by a surjective morphism $\chi \colon k[X(\G)] \to k$.

A ($k$-)\textit{character} is a $k$-point of the character variety. Any representation $\al \colon \G \to \SL_2(k)$ induces a character $\chi_\al \colon k[X(\G)] \to k$ that maps $I_\g$ to $I_\g(\al) = \Tr \al(\g)$. Such a morphism can be seen as a group homomorphism $\G \to k$, $\g \mapsto I_\g(\al)$ and we recover the standard definition of a character of the group $\G$. Finally, given any integral $k$-algebra $R$, the definition above extends to $R$-characters: an $R$-character is an $R$-point of the character variety.

The next step is to define the notions of reducible and irreducible characters. In general we are mainly interested in irreducible representations. On the other hand, it will be clear later that we cannot avoid to consider also reducible representations. 
\begin{definition}
A representation $\al \colon \G \to \SL_2(k)$ is \textit{reducible} if it preserves a one dimensional subspace in $k^2$, and \textit{irreducible} if not. More generally, given $R$ an integral $k$-algebra, a representation $\rho \colon \G \to \SL_2(R)$ is \textit{absolutely irreducible} if it is irreducible in the algebraic closure $\bar{\mK}$ of the fraction field $\mK$ of $R$.
\end{definition}
The following standard lemma allows us to define those notions directly at the level of characters in the following sense: a character is the character of an irreducible representation if and only if every representation with this character is irreducible. We will say that such a character is an irreducible character, else it is a reducible character.

\begin{lemma}\cite{CS83},\cite[Lemma 2.7]{Mar15}
\label{lem:Irred}
Let $R$ be an integral $k$-algebra (possibly $k=R$). A representation $\rho \colon \G \to \SL_2(R)$ is absolutely irreducible if and only if there exists $\g, \de \in \G$ such that $\Tr \rho(\g\de\g^{-1}\de^{-1}) \neq 2$.
\end{lemma}

We will use the notation $[\g,\de]$ for the commutator $\g\de\g^{-1}\de^{-1}$. For any $\g, \de \in \G$, we define $\Delta_{\g,\de} \in~k[X(\G)]$ as the function $I_\g^2+I_\de^2+I_{\g\de}^2-I_\g I_\de I_{\g\de} -4$. Using Lemma~\ref{lem:TraceRel}, a direct computation shows that the latter is equal to $I_{[\g,\de]}-2$. It suggests the following definition:
\begin{definition}
\label{def:Irred}
For any integral $k$-algebra $R$, an $R$-character $\chi$ is \textit{irreducible} iff there exists $\g,\de \in \G$ such that $\chi(\Delta_{\g,\de})\neq 0$. 
\end{definition}
Some reducible characters are of particular kind that we will have to exclude, namely \textit{central} characters: those are characters $\chi$ with $\chi(I_\g)^2=4$ for all $\g \in \G$.

\medskip

A consequence of Definition \ref{def:Irred} is that the reducible characters form a closed subset of the character variety $X(\G)$. Recall from Remark \ref{rem:reduced} that the character variety may have several irreducible components (it will be the case in the situation we will be interested in). Since we want to study irreducible characters, we will focus on components of the character variety that contain some of them.
\begin{definition}
\label{def:IrredType}
An irreducible component $X \subset X(\G)$ will be said \textit{of irreducible type} if it contains an irreducible character, else it is \textit{of reducible type}.
\end{definition}
Note that a component of reducible type contains exclusively reducible characters. The situation is quite different for a component of irreducible type, nevertheless irreducibility appears to be a Zarisky-open property on $X(\G)$ (see Definition~\ref{def:Irred}), hence a component of irreducible type contains an open subset of irreducible characters. In particular irreducible characters are dense in any component of irreducible type. On the other hand, a component of irreducible type may contain reducible characters. 

\medskip

Now we continue our description of characters in $X(\Gamma)$, and consider their behavior with respect to the quotient map $\pi \colon R(\G) \to X(\G)$. Clearly if two representations $\al, \al' \colon \G\to ~\SL_2(k)$ are conjugate, they define the same character $\chi_\al = \chi_{\al'}$. The converse is false in general, but true for irreducible characters, by the following proposition:
\begin{proposition}\cite[Proposition 1.5.2]{CS83}
If $\al, \al'\colon \G\to \SL_2(k)$ are representations with $\chi_\al=~\chi_{\al'}$, and if $\al$ is irreducible, then $\al$ and $\al'$ are conjugate.
\end{proposition}
In particular $\al'$ is irreducibile too. 
\begin{remark}
\label{remk:bundle}
This propostion is often summarized saying that on the irreducible part of the representation variety, the algebro-geometric quotient coincides with the topological quotient. One can be more precise saying that the map $\pi\colon R(\G)\to~X(\G)$ restricts to the irreducible part of $R(\G)$ as a principal $\operatorname{PSL}_2(k)$-bundle: it is basic linear algebra that the action by conjugation of $\SL_2(k)/\lbrace \pm I \rbrace$ is free on the set of irreducible representations.
\end{remark}
The picture is less clear in the components of reducible type. In fact, given a reducible, non-central character $\chi$ in $X(\G)$ we must distinguish two cases:

\benu
\label{prop:ReducibleChar}
\item
All representations in $\pi^{-1}\lbrace \chi \rbrace$ are abelian, in the sense that their image is an abelian subgroup of $\SL_2(k)$. In this case the character $\chi$ is said \textit{abelian}, and any two representations $\al,\al'$ with character $\chi$ are conjugate to the representation $\gamma \mapsto \bsm \la(\g)&0\\0&\la^{-1}(\g) \esm$, where $\la: \G \to k^*$ is a group homomorphism such that $\la(\g)+\la^{-1}(\g)=\chi(I_\g)$. Again, the map $\pi$ coincides here with the topological quotient, but the stabilizer $\SL_2(k)_\al$ of any representation $\al \in \pi^{-1}\lbrace \chi \rbrace$ has dimension one.
\item
There are also non-abelian representations in $\pi^{-1}\lbrace \chi\rbrace$,of the form  $\bsm \la&*\\0&\la^{-1}\esm$. In particular there are non-conjugate representations with the same character $\chi$. We say that such a character is \textit{reducible, non-abelian}. For instance reducible characters that lie in the intersection of a component of reducible type with a component of irreducible type are of this type.
\eenu

\medskip
Given an integral $k$-algebra $R$, the following theorem will allow us to lift irreducible $R$-characters to representations into $\SL_2(R)$. It is stated without proof in \cite{S93}, but a proof, written from a preceding version of \cite{S93} transmitted by Saito, can be found in the PhD thesis of the author \cite[Appendix A]{BenThesis}.
\begin{theorem}
\label{theo:Saito}
Let $R$ be an integral $k$-algebra, and let $\chi \colon B(\G) \to R$ be a morphism of $k$-algebras. Assume that $\chi(\Delta_{\g, \de})$ is invertible for some $\g, \de \in \G$, and let $A, B \in \SL_2(R)$ such that $\Tr A=\chi(I_\g)$,  $\Tr B=\chi(I_\de)$ and $\Tr AB = \chi(I_{\g\de})$. Then there exists a unique representation $\rho \colon \G \to \SL_2(R)$ whose character is $\chi$ and such that $\rho(\g)=A$ and $\rho(\de)=B$.
\end{theorem}
The following proposition is a consequence of Theorem \ref{theo:Saito}, see \cite[Proposition 3.4]{Mar15}.
\begin{proposition}
 \label{prop:Repr}
Let $\mK$ be either an algebraically closed field or a degree one extension of an algebraically closed field . The $\mK$-irreducible characters correspond bijectively to $\GL_2(\mK)$-conjugacy classes of absolutely irreducible representations $\rho\colon \G \to \SL_2(\mK)$.
\end{proposition}

\subsection{The augmented variety}
\label{sub:AugmentedVariety}
In this article we will focus on the case where $\G$ is the fundamental group of a 3-manifold $M$ with single toral boundary. We will denote the character variety $X(\pi_1(M))$ by $X(M)$, its algebra of functions $B[\pi_1(M)]$ by $B[M]$ and similarly for the boundary $\partial M$. Our main object of study will be a two-sheeted cover of the character variety, namely the augmented character variety, that we define in this subsection. The terminology (and the construction) is inspired by \cite{DG09}, it is also sometimes called the decorated character variety. This space is the space of deformations already described by Neumann--Zagier in \cite{Neumann_Zagier}. They themselves attribute its study to Thurston.

\medskip 

First, we are going to compute the character variety of the boundary of $M$.
Consider the boundary $\partial M$ of the manifold $M$, and its fundamental group $\pi_1(\partial M)$. Any character $\chi$ in $X(\partial M)$ is the character of a representation $\al \colon \pi_1(\partial M) \to~\SL_2(k)$ that can be written $\bsm \la&0\\0&\la^{-1}\esm$ up to conjugation, for some group homomorphism $\la\colon \pi_1(\partial M) \to~k^*$ such that $\la(\g)+\la^{-1}(\g) = \chi(I_\g)$. We write $\la \in H^1(\partial M, k^*)$, and $\sigma$ the involution of $H^1(\partial M, k^*)$ that turns any $\la$ into $\la^{-1}$.
\begin{remark}
Any choice of a group isomorphism $\pi_1(\partial M) \simeq \mZ^2$ induces an isomorphism $H^1(\partial M, k^*) \simeq~(k^*)^2$, and in particular it endows $H^1(\partial M, k^*)$ with a structure of affine algebraic variety. But the construction we manage to do here is intrinsic: it does not depends on this choice.
\end{remark}
\begin{proposition}
The map 
\begin{align}
\label{equa:Eigenvalue}
H^1(\partial M, k^*)/\sigma &\to X(\partial M)\\
[\la] \mapsto (\chi \colon \g &\mapsto \la(\g)+\la^{-1}(\g)) \nonumber
\end{align}
is an isomorphism of algebraic varieties. 
\end{proposition}
\begin{proof}
Consider the algebra $$C[\partial M] = k[Z_\g,  \g \in \pi_1(\partial M)]/(Z_\g Z_\de-Z_{\g\de}, Z_e-1,\g,\de\in \pi_1(\partial M)).$$ 
The elements of this algebra define functions on $H^1(\partial M, k^*)$ by $Z_\g(\la) = \la(\g)$.
We claim that this algebra is the algebra of functions of the algebraic variety $H^1(\partial M, k^*)$. To see this, observe that any identification $\pi_1(\partial M)\simeq \mZ^2$ induces the isomorphism $C[\partial M] \simeq k[X^{\pm 1},Y^{\pm1}]$, where the latter is the algebra of functions of $(k^*)^2$. 

Now let $\sigma \colon C[\partial M] \to C[\partial M]$ be defined by $\sigma(Z_\g) = Z_\g^{-1} = Z_{\g^{-1}}$. The subalgebra of invariants $C[\partial M]^\sigma$ is generated by the elements of the form $Z_\g+Z_{\g^{-1}}$, hence the morphism of $k$-algebras $B[\partial M] \to C[\partial M]^\sigma$ that maps $I_\g$ to $Z_\g + Z_{\g^{-1}}$ is an isomorphism. This proves that (\ref{equa:Eigenvalue}) is an isomorphism.
\end{proof}

\begin{definition}
\label{defi:Augmented}
The \textit{augmented character variety} is defined as the fibered product:
$$\bar{X}(M) = X(M) \times_{X(\partial M)} H^1(\partial M, k^*)$$
in other words, if $\bar{B}[M] = B[M] \otimes_{B[\pi_1(\partial M)]} C[\partial M]$, we have $\bar{X}(M) = \Spec \bar{B}[M]$.
\end{definition}
The advantage of this two-fold covering is the following: on one hand the functions of $X(M)$ are trace functions, on the other hand on $\bar{X}(M)$ we have at our disposal, for any $\g \in \pi_1(\partial M)$, two \textit{eigenvalue functions} $Z_{\g^{\pm1}}$ that map the pair $(\rho, \la)$ to an eigenvalue $\la(\g)$ of $\rho(\g)$ for any $\gamma$ in $\pi_1(\partial M)$.
\begin{remark}
\label{remk:AugmFunctions}
The algebra $\bar{B}[M]$ is generated by elements of the form $I_\g \otimes 1$ for $\g$ in $\pi_1(M)$ and $1\otimes Z_\g$ for $\g$ in $\pi_1(\partial M)$. Notice that the equality
\begin{equation}
\label{equa:AugmFunctions}
I_\g\otimes 1= 1\otimes Z_\g + 1\otimes Z_{\g^{-1}}
\end{equation}
holds for any element $\g$ in $\pi_1(\partial M)$.
\end{remark}
The following remark provides a more concrete insight on what is the augmented character variety.
\begin{remark}
Alternatively, we define the \textit{augmented representation variety} $\bar{R}(M)$ as the subvariety of $R(M) \times~H^1(\partial M,k^*)$ given by  
$$\bar{R}(M) = \lbrace (\al,\la) \in R(M) \times H^1(\partial M, k^*) | \ \la(\g)+\la(\g)^{-1} = \Tr \al(\g), \ \forall \g \in \pi_1(\partial M) \rbrace.$$ 

The group $\SL_2(k)$ acts on $\bar{R}(M)$ thus we can see the augmented character variety as the quotient $\bar{X}(M) = \bar{R}(M) /\!/ \SL_2(k)$.
\end{remark}

\subsection{The tautological representation}
\label{sub:Tautological}
In this section we define the so-called tautological representation. It has a long story in the study of character varieties, see among others \cite{CS83,CGLS,CCGLS,DFJ}, and for instance \cite{FKN} for character varieties in higher rank groups. It will be our main tool to define the Reidemeister torsion globally on the character variety.

\medskip

In this section, and from now on, we will pick $\bar{X} \subset \bar{X}(M)$ a one-dimensional component of irreducible type of the augmented character variety.
The reason why we focus on one-dimensional components is that it will allow us to define the tautological representation with entries in the function field $k(\bar{X})$, as in \cite{Mar15}. To our knowledge, it is the first occurence of such a definition, the preceding uses of tautological representations in the literature involve field extensions of $k(\bar{X})$. On the other hand it is a consequence of the work of Thurston that $\bar{X}$ is always one-dimensional if it contains the lift of a character of the holonomy representation for a hyperbolic structure on the interior of the 3-manifold M, see \cite{Thu97} and the discussion in \cite[Section 4.5]{Sha02}. Moreover, it is also proved in \cite[Section 2.3]{CCGLS} that the whole variety $X(M)$ is one-dimensional under some topological hypothesis (smallness) on the manifold $M$.

\medskip

An irreducible component $\bar{X}$ in $\bar{X}(M)$ corresponds to a minimal prime ideal $\mathfrak{p}$ of $k[\bar{X}(M)]$, in particular $k[\bar{X}]=k[\bar{X}(M)]/\mathfrak{p}$ is the algebra of function of the variety $\bar{X}$. It is an integral algebra, and the tautological morphism $$\chi_{\bar{X}}\colon k[\bar{X}(M)] \to k[\bar{X}] \to \operatorname{Frac}(k[\bar{X}]) = k(\bar{X})$$ can be seen as a $k(\bar{X})$-character (it is the generic point in the language of algebraic geometry).
\begin{lemma}
The tautological $k(\bar{X})$-character $\chi_{\bar{X}}$ is irreducible.
\end{lemma}
\begin{proof}
The tautological character specializes at any irreducible $k$-character $\chi \in \bar{X}$ as $\chi$ itself, in particular there are elements $\g,\de \in \pi_1(M)$ such that $\chi(\Delta_{\g,\de}) \neq 0$ in $k$, hence $\chi_{\bar{X}}(\Delta_{\g,\de}) \neq 0$ in $k(\bar{X})$.
\end{proof}
This tautological lemma, and the fact that $k(\bar X)$ has transcendence degree 1 over $k$ imply that we can use Proposition \ref{prop:Repr}, and we obtain the following proposition.
\begin{proposition}
\label{prop:TautRepr}
Let $\bar{X}$ be a one-dimensional irreducible component of irreducible type of $\bar{X}(M)$, and let $\chi_{\bar{X}}$ be the tautological character. There is an absolutely irreducible representation $\rho_{\bar{X}} \colon \pi_1(M) \to~\SL_2(k(\bar{X}))$, defined up to conjugation, whose character is $\chi_{\bar{X}}$. Moreover there is a tautological eigenvalue $\la_{\bar X}\colon \pi_1(\partial M) \to~k(\bar{X})^*$, defined up to inversion, such that for any $\g \in \pi_1(\partial M)$, $\lambda_{\bar{X}}(\g) = \chi_{\bar{X}}(1\otimes Z_\g)$. In particular the restricted representation $\rho_{\partial M, \bar{X}} \colon \pi_1(\partial M) \to \SL_2(k(\bar{X}))$ is diagonalizable.
\end{proposition}
\begin{proof}
The $k(\bar{X})$-character $\chi_{\bar{X}}$ restricts on $X$ to an irreducible $k(X)$-character $\chi_X \colon k[X(M)] \to ~k(X)$. Proposition \ref{prop:Repr} proves the existence of a tautological representation $\rho_X \colon \pi_1(M) \to ~\SL_2(k(X))$, defined up to $\GL_2(k(X))$-conjugation. In particular, the latter can be seen as a representation $\rho_{\bar{X}}\colon \pi_1(M) \to \SL_2(k(\bar{X}))$ defined up to conjugation by the group $\GL_2(k(\bar{X}))$.

Next, consider the group homomorphism $\la_{\bar{R}} \colon \pi_1(\partial M) \to k(\bar{R}(M))^*$ given by 
$$\la_{\bar{R}}(\g) \colon \bar{R}(M) \to k$$
$$(\rho,\la) \mapsto \la(\g)$$
Since $\la_{\bar{R}}(\g)$ is trivially $\SL_2(k)$-invariant, it induces a function on the quotient $\bar{X}(M)$, that we denote by $\la_{\bar{X}}$ after restriction to the component $\bar{X}$. In particular we have the equality $\la_{\bar{X}}(\g) +~\la_{\bar{X}}(\g^{-1}) = \chi_{\bar{X}}(I_\g\otimes 1)$. Hence up to inversion, equation (\ref{equa:AugmFunctions}) imposes $\la_{\bar{X}}(\g)$ to be equal to $\chi(1 \otimes Z_\g)$.
The last statement follows.
\end{proof}
We will frequently omit the subscript in $\rho_{\bar{X}}$ and denote the tautological representation simply by $\rho$ when the component $\bar{X}$ will be fixed. 
\begin{remark}
An important point is that the tautological representation is only defined up to conjugation. We will often write \textit{the} tautological representation, either when a representative will already be chosen or when this choice makes no difference. On the other hand, many constructions along this paper will crucially depend on the choice. For instance we will prove later that in some case, there is a subring $\cO_v$ of $k[\bar X(M)]$ and \textit{a} representative $\rho$ of the tautological representation such that $\rho \colon \pi_1(M) \to \SL_2(\cO_v)$.
\end{remark}

\subsection{Examples}
\label{sub:Examples}
In this section we write down explicit computation of (augmented) character varieties as well as tautological representations 
for the fundamental groups of the trefoil knot and of the figure-eight knot.

%
%

\subsubsection{The trefoil knot}
Let $M$ be the exterior of the trefoil knot in $\mS^3$. We are going to describe the irreducible component of irreducible type of the character variety $X(M)$. It is well-know that its fundamental group admits the presentation $\pi_1(M) = \langle a,b \vert a^2 = b^3 \rangle$, with center is isomorphic to $\mZ$ and generated by $z=a^2$. For any $\al \colon \pi_1(M) \to \SL_2(k)$ irreducible, one can show that the image of $\alpha(z)$ is central in $\SL_2(k)$, hence $\al(z) = \pm I$. But if $\al(z) = I$ then $\al(a) = -I$ and $\al$ is abelian, a contradiction. Hence $\al(z) = -I$. 

One can then assume (up to conjugation) that $\al(b) = \bsm -j&0\\0&-j^2 \esm$; where $j$ is a non trivial third root of 1. Again, because $\al$ is irreducible the right-upper entry of $\al(a)$ is not zero. Since conjugation by diagonal matrices stabilizes $\al(b)$, one can assume that $\al(a) = \bsm *&1\\ *&*\esm$. By the Cayley-Hamilton Theorem, $\al(a)^2=-I$ implies $\Tr(\al(a))=0$, and finally we have $\al(a) =~\bsm t&1\\-(t^2+1)&-t \esm$ for some $t\in k$. 

We conclude that there is a unique component of irreducible type $X$ in $X(M)$ which is isomorphic to $k$. The parameter $t$ can be seen as the function $I_{ab^{-1}}/(j-j^2)$, and the tautological representation is given by $\al \colon \pi_1(M) \to \SL_2(k(t))$. Moreover, the element $ab^{-1}$ lies in $\pi_1(\partial M)$, hence the two fold covering $\bar{X} \to X$ given by $u \mapsto \frac{u+u^{-1}}{j-j^2}=t$ is the augmented character variety. This covering ramifies twice (when $t^2(j-j^2)^2=4$), and $\bar{X}$ is isomorphic to $k^*$.

\subsubsection{The figure-eight knot}
Let $M$ be the exterior of the figure-eight knot in $\mS^3$, with fundamental group $\pi_1(M) = \langle u,v \vert uw=wv \rangle$ where $w = vu^{-1}v^{-1}u$. The first observation is that the trace functions $I_u$ and $I_v$ are equal since $u$ and $v$ are conjugated in the group, we denote this coordinate by $x$, and by $y$ we denote the function $I_{uv}$. Expanding the relation $uwv^{-1}=w$ with Lemma \ref{lem:TraceRel}, one obtains that the function ring of the character variety of the figure-eight knot is $k[X(M)] = k[x,y]/(P)$ where $P(x,y) = (x^2-y-2)(2x^2+y^2-x^2y-y-1)$. The first factor can be seen to correspond to the component of reducible type by computing $\Delta_{u,v}$, thus we denote by $X = \lbrace (x,y) \in k^2 \vert \ 2x^2+y^2-x^2y-y-1=0 \rbrace$ the plane curve defined by the second factor. It is smooth and has genus 1 by Pl\"ucker formula. Denote by $\bar{X} \to X$ the two-fold covering defined by $t+t^{-1} = x$. It ramifies four times at $\lbrace (x,y) \in X \vert \  x^2=4, y^2-5y+7=0 \rbrace$. By the Riemann-Hurwitz formula, one gets that $\bar{X}$ has genus 3.
A tautological representation is 
\begin{align*}
\rho \colon \pi_1(M) &\to \SL_2(k(\bar{X}))\\
u &\mapsto \bm t&1\\0&t^{-1} \ema\\
v &\mapsto \bm t&0\\y-t^2-t^{-2}& t^{-1} \ema
\end{align*}

\section{Differential forms, twisted cohomology and Culler-Shalen theory}
\label{sec:Twisted}
In this section we gather some facts about character varieties: in Subsection \ref{sub:DiffForm} we define the rational differential forms on $X(M)$, and rely it with the twisted homology of $M$, in Subsection \ref{sub:CS} we introduce the basic of the Culler-Shalen theory, and in Subsection \ref{sub:Twist} we perform computations in (co)-homology that will be extensively used in the upcoming sections.

\subsection{Rational differential forms on character varieties} 
\label{sub:DiffForm}
In this section we identify the first $\rho_{\bar{X}}$-twisted homology group of $M$ with the space of rational differential forms on $\bar{X}$. This result can be seen as dual of a well-know theorem due to Weil, and will be used to define the Reidemeister torsion as a rational volume form. The statement of Proposition~\ref{prop:DiffH1} and its proof are adapted from \cite{Mar15} to augmented character varieties.
\begin{notation}
The basic references for twisted homology and cohomology include \cite{Bro, DK}. For $\al \colon \pi_1(M) \to \SL_2(k)$ we will denote by $H_*(M, \Ada)$ the homology groups with coefficients in $\slf_2(k)$ where the action of the fundamental group $\pi_1(M)$ is given by 
$$\Ada \colon \pi_1(M) \xrightarrow{\al} \SL_2(k) \xrightarrow{\Ad} \Aut(\slf_2(k)).$$
Similarly, for any field extension $\mK$ of $k$ and a representation $\rho\colon \pi_1(M) \to \SL_2(\mK)$, we will denote the homology groups by $H_*(M, \Adr)$.

Finally, since we deal with Eilenberg-MacLane spaces, we will abusively use the same notation for the homology of $M$ (or of $\partial M$, or of a surface $\Sigma \subset M$) and of $\pi_1(M)$ (respectively of $\pi_1(\partial M), \pi_1(\Sigma)$).
\end{notation}

Given a representation $\al \colon \pi_1(M) \to \SL_2(k)$, a 1-cocycle $\zeta$ in $Z^1(\pi_1(M), \Ada)$ is a map $\zeta \colon \pi_1(M) \to \slf_2(k)$ satisfying the equation
\begin{equation}
\label{equa:cocycle}
\zeta(\g\de) = \zeta(\g)+\Ada(\g) \zeta(\de)
\end{equation}
It is not difficult to see that any first order deformation $\al_t$ can be written $\al + t\zeta\al$, where the map $\zeta$ satisfies (\ref{equa:cocycle}). A more precise statement is the following theorem, see for instance \cite{QL06} for a definition of Zariski tangent space.
\begin{theorem}\cite{Weil64,LM85}
\label{thm:Weil}
Let $\al \colon \pi_1(M) \to \SL_2(k)$ be an irreducible representation, such that $\chi_\al$ is a reduced (in the sense of schemes) in $X(M)$. Then there is an isomorphism $$T^{\mathrm{Zar}}_{\chi_{\al}} X(M) \simeq H^1(M, \Ada).$$
\end{theorem}

We introduce the following definition:
\begin{definition}
\label{defi:DifferentialForms}
Given a ring $A$, and a $A$-algebra $B$, we define the $B$-module of $A$-derivations $\Omega^1_{B/A}$ to be the free $B$-module generated by formal symbols  $db$, divided by the relations $\lbrace \forall a \in A, da=0$, $\forall b_1, b_2 \in B, d(b_1+b_2)=db_1+db_2$ and $d(b_1b_2)=b_1db_2+b_2db_1\rbrace$.
\end{definition}
If $X$ is an irreducible algebraic variety with function field $k(X)$, the $k(X)$-vector space $\Omega^1_{k(X)/k}$ is called the space of \textit{rational differential forms} over $X$. It is a classical fact (see \cite[Chapter~6]{QL06}) that its dimension as a $k(X)$-vector space is the dimension of $X$ as a variety over $k$.

We prove the following proposition:
\begin{proposition}
\label{prop:DiffH1}
Let $\bar{X}\subset \bar{X}(M)$ be a one-dimensional component of irreducible type of the augmented character variety, with function ring $k[\bar{X}] = k[\bar{X}(M)]/\mfp$, and we fix a tautological representation $\rho\colon \pi_1(M) \to ~\SL_2(k(\bar{X}))$ such that $\rho(\partial M)$ is diagonal. There is an exact sequence of $k(\bar X)$-vector spaces
\begin{equation}
\label{equa:exactDiff}
\mfp/\mfp^2 \otimes_{k[\bar{X}]} k(\bar{X}) \to H_1(M, \Adr) \to \Omega^1_{k(\bar{X})/k} \to 0
\end{equation}
\end{proposition}
\begin{proof}
The proof of the proposition will follow from classical arguments of algebraic geometry once we have proved the following claim:
\begin{claim}
There is an isomorphism 
\begin{equation*}
\Omega^1_{\bar{B}[M]/k} \otimes_{\bar{B}[M]} k(\bar{X}) \simeq H_1(M, \Adr).
\end{equation*}
\end{claim}
\begin{proof}[Proof of the claim]
For any $\g$ in $\pi_1(M)$, we denote by $\rho(\g)_0 \in \slf_2(k(\bar{X}))$ the trace-free matrix given by $\rho(\g)-~\frac{1}{2}\Tr(\rho(\g))I$. The space $C_1(M, \Adr)$ is generated by elements of the form $\xi\otimes[\gamma]$, with $\xi \in \slf_2(k(\bar{X})), \g \in \pi_1(M)$. 

We construct a morphism of $\bar{B}[M]$-modules $$\Omega^1_{\bar{B}[M]/k} \to H_1(M,\Ad\circ\rho)$$
$$d(I_\g \otimes 1) \mapsto \rho(\g)_0 \otimes [\g]$$
$$d(1\otimes Z_\g) \mapsto \bsm \frac \la 2 & 0\\0 & -\frac \la 2 \esm \otimes [\g] $$
where $\lambda, \lambda^{-1}$ are the eigenvalues of the matrix $\rho(\g)$.

Using $\partial \xi \otimes [\gamma] = \rho(\gamma)^{-1} \xi \rho(\gamma) - \xi$, one gets that $d(I_\gamma \otimes 1)$ and $d(1\otimes Z_{\gamma})$ are well-defined cocycles (remember that $\rho(\gamma)$ is diagonal for $\gamma \in \pi_1(\partial M)$).
Moreover, using the formula $\partial \xi \otimes [\gamma,\delta] = \rho(\gamma)^{-1}\xi \rho(\gamma) \otimes[\delta] -\xi \otimes [\gamma\delta] +\xi \otimes [\gamma]$ one can show that $d(I_\gamma I_\delta \otimes 1)$ and $d(I_{\gamma\delta} \otimes 1 + I_{\gamma\delta^{-1}} \otimes 1)$ are mapped to the same element in $H_1(M, \Adr)$.

In addition, it induces
$$d(1\otimes Z_{\g^{-1}}) \mapsto \bsm \frac {\la^{-1}} {2} & 0\\0 & -\frac {\la^{-1}}{2} \esm \otimes [\g^{-1}]=  \bsm -\frac {\la^{-1}} {2} & 0\\0 & \frac {\la^{-1}}{2} \esm \otimes [\g],$$
hence this morphism is well-defined.
It induces a $k(\bar{X})$-linear map $$\Psi\colon~ \Omega^1_{\bar{B}[M]/k} \otimes k(\bar{X}) \to H_1(M,\Ad\circ\rho).$$

To construct the reciprocal morphism, we define $\bar{\Lambda} = k(\bar{X}) \oplus \varepsilon \Omega^1_{\bar{B}[M]/k} \otimes k(\bar{X})$, and following \cite{Mar15} we define the map $\varphi: \bar{B}[M] \to \bar{\Lambda}$ given by 
$$I_\g\otimes 1 \mapsto I_\g \otimes 1 + \varepsilon d(I_\g \otimes 1)$$
$$1\otimes Z_\g \mapsto 1\otimes Z_\g + \varepsilon d(1\otimes Z_\g)$$
By Theorem \ref{theo:Saito}, we can produce a representation $\rho_\varepsilon: \pi_1(M) \to \SL_2(\bar{\Lambda})$ such that $\chi_{\rho_\varepsilon} = \varphi$. 

Now using the fact that the vector $\frac d {d\varepsilon}\rho_\varepsilon(\gamma) \rho(\gamma)^{-1}$ is trace free one can check that the map $\xi\otimes[\g] \mapsto \frac d {d\varepsilon} \Tr(\xi\rho_\varepsilon(\g)\rho(\g)^{-1})$ is a left-section of the morphism $\Psi$ above. Indeed, it is also a right-section since $H_1(M, \Adr)$ is linearly generated by cycles of the form $\rho(\gamma)_0 \otimes [\gamma]$.
\end{proof}
Now we deduce the proposition from the claim. First, we deduce from \cite[Chapter 6, Proposition 1.8, (c)]{QL06} the isomorphism $\Omega^1_{k[\bar{X}]/k}\otimes ~k(\bar{X})\simeq \Omega^1_{k(\bar{X})/k}$, and then we know from \cite[Chapter 6, Proposition 1.8, (d)]{QL06} that the map $$\Omega^1_{\bar{B}[M]/k} \otimes k[\bar{X}] \to \Omega^1_{k[\bar{X}]/k}$$ is onto, with kernel $\mfp/\mfp^2$, and the proposition follows by tensoring with the field $k(\bar{X})$.
\end{proof}

\begin{corollary}
\label{cor:Cotangent}
Let $\bar{X}$ be a one-dimensional irreducible component of irreducible type of $\bar{X}(M)$ corresponding to a minimal prime $\mfp$. If $\mfp/\mfp^2 \otimes_{k[\bar{X}]}~k(\bar{X})$ is trivial, then there is an isomorphism
$$H_1(M, \Adr)\simeq \Omega^1_{k(\bar{X})/k}.$$
\end{corollary}
\begin{remark}\label{remk:EssRed}
We will say that the component $\bar{X}$ is \textit{essentially reduced in $\bar{X}(M)$}  if the hypothesis $\mfp/\mfp^2 \otimes_{k[\bar{X}]} k(\bar{X})=0$ holds. If $X(M)$ is scheme-reduced, then any component $X$ of $X(M)$ is essentially reduced. For sake of generality, in the rest of this article we will keep this minimal hypothesis of $X$ being essentially reduced, but the reader can with few losses think that $X(M)$ is supposed to be reduced. The fact is that this hypothesis is global on $X(M)$, whereas we deal only with a fixed component $X \subset X(M)$.

\end{remark}

\subsection{Smooth projective model and Culler-Shalen theory}
\label{sub:CS}
In this subsection we introduce the material from \cite{CS83} (see \cite{Til03,Sha02} for expository notes on the topic) to construct embedded essential surfaces in 3-manifolds $M$ from the action of their fundamental group $\pi_1(M)$. 

\medskip

Given a curve $X$ over an algebraically closed field $k$ of characteristic zero, there is a unique smooth projective curve $\widehat{X}$, up to isomorphism, that is birational to $X$ (see \cite{Ful} for a detailed exposition on algebraic curves). It is called the \textit{smooth projective model} of $X$.
One way to define it is to consider the set of discrete $k$-valuations on the function field $k(X)$, with the cofinite topology. A discrete $k$-valuation is a group epimorphism $v \colon k(X)^* \to \mZ$ such that $v(f+g) \ge \min(v(f),v(g))$ and $v(k^*)=0$. It is extended to $k(X)$ by $v(0)=\infty$. The birational map $\nu: \widehat{X} \to X$ can be described as follows: to any smooth point $x \in X$ corresponds a unique valuation $v_x$ on $k(X)$ that maps a rational function $P$ on the vanishing order of $P$ at $x$. If $x$ is not smooth, it may exist several ways to define $v_x$ as such. On the other hand there exist valuations $v$ on $k(X)$ that do not correspond to a point $x$ of $X$, such valuations are called \textit{ideal} valuations, or ideal points of $\widehat{X}$. Ideal points will play a crucial role in Culler-Shalen theory.

We will denote by $\cO_v = \lbrace P \in k(X) \vert \ v(P) \ge 0 \rbrace$ the \textit{valuation ring} of $v$. It has many advantageous properties, for instance it is a principal ring with a unique maximal ideal that we will denote by $(t)$. Such a choice of $t$ will be called an \textit{uniformizing element}, it is characterized by the fact that $v(t)=1$. Moreover, any ideal of $\cO_v$ is of the form $(t^n)$, for some natural number $n$. To give an insight on what the ring $\cO_v$ looks like, one can think about the ring of formal series $k[[t]]$; in fact, the valuation rings we encounter in the theory have $k[[t]]$ as a natural completion. Finally, the \textit{residual field} is the field $\cO_v/(t)$, in our context it is isomorphic to $k$. Note that $v$ is ideal iff $k[X]$ is not contained in $\cO_v$.

\medskip

\begin{notation}
In the sequel of this article, we will work with the smooth projective model of a one-dimensional component of irreducible type $\bar{X}$ essentially reduced in $\bar{X}(M)$, and to avoid too many superscripts, we will denote it by $Y$ rather than $\widehat{\bar{X}}$ that would be quite unaesthetic.
\end{notation}

\begin{remark}
The morphism $\nu$ induces a field isomorphism $k(Y) \simeq k(\bar{X})$, in particular everything we proved, up to now, concerning $k(\bar{X})$ remains true replacing it by $k(Y)$.
\end{remark}

To any valuation $v$ on $k(\bar{X})$ (from now on we will write "to any point $v$ in $Y$"), the Culler-Shalen theory associates an action of the fundamental group $\pi_1(M)$ on a tree $T_v$. The vertices of $T_v$ are in bijection with the homothety classes of $\cO_v$-lattices $L \subset k(Y)^2$, and the action is given by the tautological representation $\rho \colon \pi_1(M) \to \SL_2(k(Y))$. 

The tree $T_v$ associated to the valued field $(k(Y),v)$ is called the \textit{Bass-Serre tree}. We recommend the foundational reference \cite{Ser77} as well as Chapter 3 of \cite{Sha02} for a detailed treatment of this theory. Since for any (class of) lattice $L$ in $k(Y)^2$, one can fix a basis of $k(Y)^2$ such that $L \simeq \cO_v^2$, one deduces the following lemma from the fact that the stabilizer of $\cO_v^2$ in $\SL_2(k(Y))$ is precisely $\SL_2(\cO_v)$:

\begin{lemma}
The image of the representation $\rho\colon \pi_1(M) \to \SL_2(k(Y))$ fixes a vertex of the tree $T_v$ if and only if it is conjugated to a subgroup of $\SL_2(\cO_v)$
\end{lemma}

In this case, if we chose the tautological representation $\rho \colon \pi_1(M)\to \SL_2(\cO_v)$, such a representative will be said \textit{convergent} at $v$, since for any $\g$ in $\pi_1(M)$, the evaluation at the point $v$ of the entries of $\rho(\g)$ is finite. The following proposition characterizes valuations $v \in Y$ such that there exists a tautological representation that converges at $v$.
\begin{proposition}
\label{prop:FixIdeal}
There is a convergent tautological representation at $v$ if and only if $v$ is not an ideal point.
\end{proposition}
\begin{proof}
The point $v$ being finite is equivalent to the fact that for any $\g$ in $\pi_1(M)$, the trace $\Tr(\rho(\g))$ lies in $\cO_v$ since $k[\bar{X}] \subset \cO_v$. We claim that it is equivalent to $\rho(\g)$ to be conjugated to a matrix in $\SL_2(\cO_v)$: it is clear if $\rho(\g) =\pm I$, else there exists a vector $V \in k(Y)^2$ such that $\lbrace V, \rho(\g) V \rbrace$ is a basis of $k(Y)^2$, in this basis $\rho(\g)$ is $\bsm 0&1\\-1&\Tr\rho(\g)\esm$ and the claim follows. Now the proof of Proposition \ref{prop:FixIdeal} is an immediate consequence of the following lemma.
\end{proof}
\begin{lemma}(See \cite[Corollaire 3 p.90]{Ser77}, \cite[Lemma 1.3.7]{BenThesis}.)
If $G$ is a subgroup of $\SL_2(k(Y))$ such that any element $g \in G$ fixes a vertex of the Bass-Serre tree $T_v$, then the whole group $G$ fixes a vertex of $T_v$.
\end{lemma}

Let $v \in Y$ a valuation on $k(Y)$, and $T_v$ the Bass-Serre tree associated to the valued field $(k(Y),v)$. In \cite{CS83}, M. Culler and P. Shalen construct a $\pi_1(M)$-equivariant simplicial map $f \colon \widetilde{M} \to T_v$. Denoting by $E$ the set of mid-edges in $T_v$, they claim that $f^{-1}(E)$ is an invariant surface in $\widetilde{M}$ (possibly empty). In particular it defines a surface $\Si_v \subset M$ in the quotient, which is said to be \textit{dual} to the action.
\begin{definition}
A surface $\Si$ in a 3-manifold $M$ is said \textit{essential} if
\benu
\item
The surface $\Si$ is properly embedded, that is $(\Si, \partial \Si) \subset (M, \partial M)$ is an embedding.
\item
$\Si$ is oriented.
\item
No component of $\Si$ is a two-sphere, or is parallel to the boundary.
\item
For each component $\Si_i$ of $\Si$, the induced homomorphism $\pi_1(\Si_i) \to \pi_1(M)$ is one-one.
\eenu
\end{definition}

Finding essential surfaces in 3-manifolds is a deeply-studied question. For instance, for $M$ a 3-manifold with toral boundary, knowing if there exists a separating essential surface $\Si \subset M$ with non empty boundary was known as the \textit{weak Neuwirth conjecture}.

In the Culler-Shalen construction above, the dual surface $\Si_v$ can be rendered essential in the manifold $M$. Among many consequences, it provided a proof of the weak Neuwirth conjecture (\cite{CS84,Sha02}).

\begin{remark}
The construction above depends on many choices, in particular a dual surface is not unique. In this paper we will work with a fixed dual surface satisfying certain conditions, and we will obtain certain inequalities involving the Euler characteristic of this surface. If several surfaces would satisfy the conditions, then we might pick the one for which the result is optimal.
\end{remark}

The following proposition emphasizes the importance of ideal points in the Culler-Shalen construction; it follows easily from Proposition \ref{prop:FixIdeal}.
\begin{proposition}
\label{prop:Nonempty}
If $v \in Y$ is an ideal point, then there exists a non-emtpy dual surface $\Si_v$.
\end{proposition}

In what follows we will focus on the following situation (corresponding essentially to the Neuwirth conjecture): $v \in Y$ is an ideal point such that the dual surface $\Si_v$ is essential, separating, $\partial \Si_v \neq \emptyset$. We furthermore assume that the surface $\Si_v$ consists of $n$ parallel copies $\Si_i$ with $M \setminus \Si_i$ homeomorphic to the union of two handlebodies $M_1$ and $M_2$ (in this case the surface $\Si_v$ is said to be \textit{free}). It is described in the following picture.
\begin{figure}[h]
\begin{center}
\def\svgwidth{0.6\columnwidth}
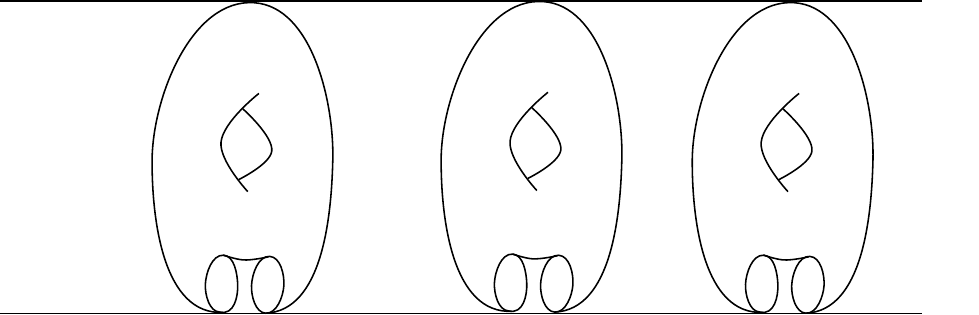
\caption{\label{fig:split surface} The splitting $M = M_1 \cup_{V(\Si)} M_2$. Note that $\partial V(\Sigma) = \Sigma_1 \cup \Sigma_n$. The quotient of the corresponding Bass-Serre tree is drawn.}
\end{center}
\end{figure}
Hence $V(\Si)\simeq \Si_1 \times [0,1]$ is a thickening of any $\Si_i$. The assumption of $\Si_v$ consisting of parallel copies is slightly more general than $\Si_v$ being connected: in general there is no reason for $\Si_v$ to be connected, but it is known that the boundary curves $\partial \Si_v$ are parallel in the peripheral torus $\partial M$. 

\medskip

We fix a basepoint $p \in \Si_1$, and we denote by $\pi_1(M_1)$, $\pi_1(M_2)$ and $\pi_1(\Si)$ respectively the fundamental groups of $M_1$, $M_2 \cup V(\Si)$ and of $\Si_1$ relatively to the point $p$. Denote by $i_1, i_2$ the maps induced by inclusion $\pi_1(M_1) \to \pi_1(M)$. The Seifert-Van Kampen Theorem induces the amalgamated product
\begin{equation}
\label{equa:Amalgam}
\pi_1(M) = \pi_1(M_1) \ast_{\pi_1(\Si)} \pi_1(M_2)
\end{equation}
and $i_\Si$ is the composition $\pi_1(\Si) \to \pi_1(M_1) \xrightarrow{i_1} \pi_1(M)$ or equivalently $\pi_1(\Si) \to \pi_1(M_2) \xrightarrow{i_2}\pi_1(M)$.

\begin{notation}
For $\rho \colon \pi_1(M) \to \SL_2(k(Y))$, one denotes by $\rho_1$ the composition  $$\rho_1\colon \pi_1(M_1) \xrightarrow{i_1}\pi_1(M) \xrightarrow{\rho} \SL_2(k(Y)),$$ by $\rho_2 \colon \pi_1(M_2) \xrightarrow{i_2} \pi_1(M) \xrightarrow{\rho} \SL_2(k(Y))$ and by $\rho_\Si \colon \pi_1(\Si) \xrightarrow{i_\Si} \pi_1(M) \xrightarrow{\rho} \SL_2(k(Y))$.
\end{notation}
Since $v$ is an ideal point, the tautological representation $\rho$ cannot be convergent, on the other hand, the surface $\Si$ is precisely the locus in $M$ where it diverges, in the sense that $\rho$ can be chosen to converge on each pieces $\pi_1(\Si)$, $\pi_1(M_i), i=1,2$. A more precise statement is given in the following lemma.
\begin{lemma}
\label{lem:Converge}
Let $v \in Y$ be an ideal point, there is a tautological representation $\rho\colon \pi_1(M) \to \SL_2(k(Y))$ such that $\rho_1$ is convergent and that $\rho_\Sigma$ is residually reducible. Moreover, there is a convergent representation $\rho_2'\colon \pi_1(M_2)\to \SL_2(\cO_v)$ such that $\rho_2 = U_n \rho_2' U_n^{-1}$, with $U_n = \bsm t^n & 0\\0 & 1 \esm$.
\end{lemma}
\begin{proof}
Let $s_1 \in T_v$ be a vertex in the Bass-Serre tree that is fixed by $\pi_1(M_1)$, and fix a basis such that it corresponds to the lattice $\cO_v^2$. We claim that there is a vertex $s_2 \in T_v$, fixed by $\pi_1(M_2)$, such that $d(s_1,s_2)=n$. Moreover,  we may assume that in this basis $s_2$ has a representative of the form $t^n\cO_v \oplus \cO_v$. Let us prove the claim.

The path linking $M_1$ to $M_2$ in Figure~\ref{fig:split surface}, lifted in the universal cover, is by construction (see \cite[Section 2]{Sha02}) the pull-back of a segment of length $n$ linking the vertex $s_1$ to a vertex $s_2$ in the Bass-Serre tree, such that $\pi_1(M_2)$ stabilizes $s_2$. We prove inductively that there is a basis of $k(Y)^2$ such that $s_2$ represents the lattice $t^n \cO_V \oplus \cO_v$. To do so, denote by $s$ the first vertex on the segment from $s_1$ to $s_2$, we prove that $s$ represents $t\cO_V \oplus \cO_v$ in an appropriate basis.

Since $d(s_1, s)=1$, any lattice $L$ in the class of $s$ can be written $a\cO_v \oplus b \cO_v$, with $a,b$ in $\cO_v$, $|v(b)-v(a)|=1$. Up permute the basis elements, one can suppose that $v(b)-v(a)=1$. Now in the homothety class of $s$ we pick $L$ such that $v(b)=1$ and $v(a)=0$. Up to change the basis by a matrix in $\SL_2(\cO_v)$, one can take $a=1$, $b=t$, so that $s_1$ corresponds to $\cO_v \oplus t \cO_v$. Finally, using this argument inductively one can fix a basis such that $s_2$ is the lattice $\cO_v \oplus t^ncO_v$.

The first observation is that $\rho_1(\pi_1(M_1)) \subset \SL_2(\cO_v)$ because it stabilizes $\cO_v^2$. Since $\rho_\Sigma$ fixes the first edge of the segment $[s_1s_2]$, in this basis it fixes the lattices $\cO_v^2$ and $t\cO_v \oplus \cO_v$, hence for all $\g \in \pi_1(\Sigma)$, $\rho_\Sigma(\g) = \bsm a(\g)& b(\g)\\c(\g) & d(\g) \esm$, with $c(\g) \in (t)$, hence $\brho_\Sigma$ is reducible. 

Let $\rho'_2 = U_n^{-1} \rho_2 U_n$, then $\rho'_2 \cdot s_1  = U_n^{-1} \rho_2 \cdot s_2 = U_n^{-1} \cdot s_2 = s_1$ and we have proved that the representation $\rho'_2$ converges.
\end{proof}

\subsection{Some computations of twisted cohomology groups}
\label{sub:Twist}
In this subsection we collect some technical lemmas that will be used along this article on the $\rho$-twisted cohomology of $M$.

\medskip

First we compute the cohomology of $M$ with coefficients in the function field $k(Y)$. The matrix $\bsm 1&0\\0&-1\esm$ of $\slf_2(k(Y))$ will be denoted by $H$. The group homomorphism $\pi_1(\partial M) \to \pi_1(M)$ induces a morphism $r^*\colon H^*(M, \Adr) \to H^*(\partial M, \Adr)$.
\begin{proposition}
\label{prop:CohomDimension}
For $i=0$ or $i \ge 3$, the $k(Y)$-vector spaces $H^i(M, \Adr)$ are trivial, and there is a natural isomorphism $H^2(M, \Adr) \simeq k(Y)$.
\end{proposition}

\begin{proof}
Recall that the tautological representation $\rho$ is irreducible. In particular the vector space $H^0(M, \Adr)$, which is by definition the space of $\Adr$-invariant matrices in $\slf_2(k(Y))$, is trivial. It is well-know that a connected 3-manifold with non-empty boundary has the same homotopy type as a 2-dimensional CW complex, in particular it has no homology in dimension 3 or higher, and the first statement is proved.

Since $\partial M$ is a torus one has $\chi(M) = \frac 1 2 \chi(\partial M) = 0$, hence $$\dim(H^1(M, \Adr))  = \dim (H^2(M,\Adr)).$$ From Corollary \ref{cor:Cotangent} together with the fact at the end of Definition \ref{defi:DifferentialForms}, this dimension is known to be the dimension of $Y$ as a variety over $k$, hence $\dim(H^2(M, \Adr))=1$.

Now consider the following part of the long exact sequence of the pair $(M,\partial M)$ in twisted cohomology
\begin{equation}
\label{equa:LongSequ}
H^2(M, \Adr) \xrightarrow{r^*} H^2(\partial M, \Adr) \to H^3(M, \partial M, \Adr) 
\end{equation}
By Poincar\'e-Lefschetz duality, the right-hand side term in the sequence (\ref{equa:LongSequ}) is isomorphic to $H_0(M, \Adr)$ that is trivial, hence the map $r^*$ is onto. Now we have $\dim(H^2(M, \Adr))=1$, so that it is enough to show that $H^2(\partial M, \Adr)$ is non-trivial, it will follow that the map $r^*$ is an isomorphism of $k(Y)$ vector spaces. This last claim is true because of the duality $$H^2(\partial M, \Adr) \simeq H_0(\partial M, \Adr) \simeq H^0(\partial M, \Adr)^*$$ (here $H^0(\partial M, \Adr)$ is non-trivial because $\pi_1(\partial M)$ is abelian). Hence $r^*$ is an isomorphism.

In particular it shows that $H^0(\partial M, \Adr)$ has dimension $1$. By Proposition \ref{prop:TautRepr} one can assume that $\rho(\pi_1(\partial M))$ is a diagonal subgroup of $\SL_2(k(Y))$, hence $H^0(\partial M, \Adr)$ is spanned by the diagonal matrix $H$. Finally the natural homomorphism 
\begin{align}
\label{equa:NaturalIsom}
H^2(M, &\Adr) \to k(Y)\\ \nonumber
\eta &\mapsto \Tr(\eta[\partial M] H)
\end{align}
is an isomorphism: it is the composition of the isomorphism $r^*\colon H^2(M, \Adr) \simeq H^2( \partial M, \Adr)$ above with the identification $H^2(\partial M, \Adr) \simeq k(Y)$ induced by Poincar\'e duality.
\end{proof}

In the second part of this subsection, we focus on the twisted cohomology with coefficients in  $\slf_2(\cO_v)$. We fix $v \in Y$ a finite point, by Proposition \ref{prop:FixIdeal} one can fix a convergent tautological representation $\rho\colon\pi_1(M) \to \SL_2(\cO_v)$, and define the complex of $\cO_v$-modules $C^*(M, \Adr)_v$ of twisted cohomology with coefficients in $\slf_2(\cO_v)$.

We define the \textit{residual representation} $\brho \colon \pi_1(M) \xrightarrow{\rho} \SL_2(\cO_v) \xrightarrow{mod (t)} \SL_2(k)$. We will say that a convergent tautological representation is \textit{residually reducible} (respectively \textit{abelian, central}) if the residual representation is reducible (resp. abelian, central). The residual complex is the complex of $k$-vector spaces $C^*(M, \Adbr)$, and it will be used extensively in the sequel.

The Reidemeister torsion will be seen to be related with the torsion part in the $\cO_v$-modules $H_*(M, \Adr)_v$. We will use the following lemma:
\begin{lemma}
\label{lem:H0}
If the tautological representation $\rho$ is not residually central at a finite point $v$ in $Y$, then it can be chosen such that the module $H_0(M, \Adr)_v$ is trivial.
\end{lemma}
\begin{remark}
\label{remk:Invariance}
Note that in the statement of Lemma \ref{lem:H0}, it is implicit that being residually central for the tautological representation  does not depend on the choice of a convergent representative $\rho\colon \pi_1(M) \to \SL_2(\cO_v)$, because being residually central is a property of the tautological character. That is not true for the property of being residually abelian, in particular this lemma shows that the $\cO_v$-modules $H_i(M, \Adr)_v$ depend on the choice of $\rho$ in its $\GL_2(k(Y))$-conjugacy class. The explanation why at the end, the computation of the torsion, that involves those $\cO_v$-modules, will not, is because the torsion is an invariant of the complex of $k(Y)$-vector spaces $C^*(M, \Adr)$, which does not depend on the $\GL_2(k(Y))$-conjugacy class of the representation~$\rho$.
\end{remark}

Before giving a proof of Lemma \ref{lem:H0}, we discuss briefly the hypothesis that the tautological representation is not residually central. For our purpose, it will follow from the hypothesis that the manifold $M$ has the rational homology of a circle, as showed now:
\begin{lemma}
\label{lem:NonCentral}
Let $M$ a 3-manifold with toral boundary and first Betti number equal to 1. If $X$ is a component of irreducible type of the character variety $X(M)$, then it does not contain any central character.
\end{lemma}
\begin{proof}
Assume $\chi$ is central, then any representation $\al \colon \pi_1(M) \to \SL_2(k)$ whose character is $\chi$ is conjugated to a representation of the form $\al(\g) = \pm \bsm 1 & \varphi(\g)\\0 & 1 \esm$, for some $\varphi \colon \pi_1(M) \to \mZ$. Since $b_1(M) =1$, every such non-trivial representations are conjugated, and the dimensional argument of \cite[Lemme 3.9, (iii)]{Porti97} can be applied: the fiber $\pi^{-1}(\lbrace \chi \rbrace)$ has dimension two, what contradicts the fact that $\chi$ lies in a component of irreducible type $X$, where the fibers of the quotient map $\pi$ have dimension at least 3.
\end{proof}
A more general statement, and a systematic study of deformation theory of central characters, can be found in \cite{MarcheWolff}.

\medskip

Lemma \ref{lem:H0} will promptly follow from the following lemma:
\begin{lemma}
\label{lem:NotAbelian}
Let $v \in Y$ be a finite valuation, and $\rho \colon \pi_1(M)\to \SL_2(\cO_v)$ a convergent tautological representation. If $\rho$ is not residually central, then it can be chosen to be not residually abelian.
\end{lemma}

\begin{proof}
We give a tree-theoretical argument, in the spirit of this article.
Let $T_v$ be the Bass-Serre tree associated to $v$, and consider the action of the fundamental group $\pi_1(M)$ on $T_v$. Let $T'_v$ be the subtree of fixed points in $T_v$. Since the representation $\rho$ is convergent, the tree $T'_v$ is not empty.

Now we show that the tree $T'_v$ is finite: if not it would contain a half-line, hence from \cite[p.107]{Ser77} the tautological representation would fix a line in the completion $\hat{\cO}_v^2$ and it would contradict the irreducibility of $\rho$.

We claim that the tree $T'_v$ is a segment. By way  of contradiction, assume that it contains a vertex $s$ of valence at least~3:
\begin{center}
\begin{tikzpicture}
\node[circle,fill,label=$s$,inner sep=1pt] (s) at (0,0){};
\node[circle,fill,label=$t$,inner sep=1pt] (t) at (1,0){};
\node[circle,fill,label=$u$,inner sep=1pt] (u) at (-0.7,0.7){};
\node[circle,fill,label=$v$,inner sep=1pt] (v) at (-0.7,-0.7){};
\draw (u)--(s)--(t);
\draw (s)--(v);
\end{tikzpicture}
\end{center}
Then $t,u$ and $v$ represents three distincts lattices at distance 1 of $s$. They correspond to points in $k\mathbb{P}^1$ (see \cite[Section 3.8]{Sha02}) that are fixed by the residual representation $\bar{\rho}$. Hence the whole image $\bar{\rho}(\pi_1(M))$ in included in  $\lbrace \pm I \rbrace$ and it contradicts the fact that the tautological representation is not residually central.

Finally we proved that $T'$ is of the following form: 
\begin{center}
\begin{tikzpicture}
\node[circle,fill,label=$s_0$,inner sep=1pt] (s0) at (0,0){};
\node[circle,fill,label=$s_1$,inner sep=1pt] (s1) at (1,0){};
\node[circle,fill,label=$s_2$,inner sep=1pt] (s2) at (2,0){};
\node[] (s) at (3,0){...};
\node[circle,fill,label=$s_{n-1}$,inner sep=1pt] (s3) at (4,0){};
\node[circle,fill,label=$s_n$,inner sep=1pt] (sn) at (5,0){};
\draw (s0)--(s1)--(s2)--(s)--(s3)--(sn);
\end{tikzpicture}
\end{center}
Let us fix a basis such that $s_0$ represents the lattice $\cO_v^2$. 
With the same inductive argument as in the proof of Lemma \ref{lem:Converge}, we prove that there is a basis of $k(Y)^2$ such that $s_i$ represents the lattice $\cO_v \oplus t^i\cO_v$ for any $i=0, \ldots, n$.

Hence for any $\g$ in  $\pi_1(M)$, we have $\rho(\g) =\bsm a(\g) & b(\g)\\c(\g) & d(\g) \esm \in~\SL_2(\cO_v)$, with $c(\g) \in (t^n)$. We conclude by noting that for some $\g \in \pi_1(M)$, we have $b(\g) \in \cO_v^\ast$; because if not there should be an other fixed point at the left of $s_0$. Its ends up the proof of the lemma since this choice of tautological representation $\rho$ is not residually abelian.
\end{proof}

\begin{proof}[Proof of Lemma \ref{lem:H0}]
We have isomorphisms of $k$-vector spaces
\begin{equation}
\label{equa:H0}
H_0(M, \Adr)_v \otimes_{\cO_v} \cO_v/(t) \simeq H_0(M, \Adbr) \simeq H^0(M, \Adbr)
\end{equation}
where in Equation (\ref{equa:H0}) the first isomorphism comes from the Universal Coefficient Theorem (coefficient change) and the second from Universal Coefficient Theorem (duality). From Lemma \ref{lem:NotAbelian}, one can chose a non-abelian tautological representation $\rho$, thus $H^0(M, \Adbr) = \lbrace 0 \rbrace$ hence $H_0(M, \Adr)_v \otimes_{\cO_v} \cO_v/(t)$ is trivial. It implies that the $\cO_v$-module $H_0(M, \Adr)_v$ is torsion-free. Now $\rho$ is irreducible (in particular non-abelian) hence the $k(Y)$-vector space $H_0(M, \Adr)$ is trivial. Again by the UCT, $H_0(M, \Adr)_v \otimes_{\cO_v} k(Y)\simeq H_0(M, \Adr)$, hence $H_0(M, \Adr)_v$ has trivial rank, and is the trivial module.
\end{proof}

We finish this section by a computation of the $\cO_v$-modules $H_*(M, \Adr)_v$ and $H^*(M, \Adr)_v$.
\begin{lemma}
\label{lem:Homology}
Let $\rho$ be such that Lemma \ref{lem:H0} holds, there are the following isomorphisms:
\begin{align}
\label{1}
& H^0(M, \Adr)_v = \lbrace 0 \rbrace \\
\label{2}
& H_1(M, \Adr)_v \simeq \Omega^1_{\bar{B}[M]/k} \otimes \cO_v \text{ if $\rho$ is residually irreducible }\\
\label{3}
& H^1(M, \Adr)_v \simeq \cO_v\\
\label{4}
& H_2(M, \Adr)_v \simeq \cO_v\\
\label{5}
& H^2(M, \Adr)_v \simeq \cO_v \oplus T_1(M)\\
\label{6}
&H_i(M, \Adr)_v = H^i(M, \Adr)_v = \lbrace 0 \rbrace \text{ else }
\end{align}
where $T_1(M)$ denotes the torsion part of the $\cO_v$-module $H_1(M, \Adr)_v$.
\end{lemma}
\begin{proof}
Since $M$ has the same homotopy type as a two-dimensional CW-complex, (\ref{6}) holds.
For any~$i$, the rank of the modules $H_i(M, \Adr)_v$ and $H^i(M, \Adr)_v$ is determined by Proposition~\ref{prop:CohomDimension}. This rank is $0$ for $i=0$, and $1$ for $i=1,2$. Since $H^0(M, \Adr_v)$ is free, (\ref{1}) holds, and the UC exact sequence $$0 \to \Ext(H_0(M, \Adr)_v,\cO_v) \to H^1(M, \Adr)_v \to \Hom(H_1(M, \Adr)_v, \cO_v) \to 0$$ together with the fact that $H^1(M, \Adr_v)$ has rank one give (\ref{3}). The module $H_2(M, \Adr)_v$ is free because there are no 3-dimensional chains, hence (\ref{4}), and then (\ref{5}) follows from the UC exact sequence $$0 \to \Ext(H_1(M, \Adr)_v,\cO_v) \to H^2(M, \Adr)_v \to \Hom(H_2(M, \Adr)_v, \cO_v) \to 0.$$
Finally the proof of (\ref{2}) follows closely the claim in the proof of Proposition \ref{prop:DiffH1}: the same construction provides a morphism of $\cO_v$-modules $\Omega_{\bar{B}[M]/k} \otimes \cO_v \to H_1(M, \Adr)_v$. Since $\rho$ is residually irreducible, there exists $\g, \de \in \pi_1(M)$ such that $\chi_\rho(\Delta_{\g,\de}) \in \cO_v^*$ and Theorem \ref{theo:Saito} applies.
\end{proof}

\section{The Reidemeister torsion form}
\label{sec:Reidemeister}
In this section we define the main object of this article: the Reidemeister torsion form.
In Subsection \ref{sub:DefiTorsion} we give a quick overview of the Reidemeister torsion theory. In Subsection \ref{sub:TorsionForm} we explain how it gives rise to a rational differential form on the augmented character variety.
\subsection{Reidemeister torsion}
\label{sub:DefiTorsion}
In this subsection we give various definitions of the Reidemeister torsion. There are many good references on the topic, among them we follow \cite[Appendix A]{GKZ94}, \cite{Mil66} and particularly \cite{Porti97} which corresponds to our situation. The convention followed here, namely where we start to index the alternating product in the definition of the determinant of a complex, is the one of \cite{GKZ94} but is the opposite of the convention followed in \cite{Mil66, Porti97}.

\medskip

\subsubsection{Definition of the torsion}
\label{subsub:Definition}
Given a finite complex $C^*$ of $k$-vector spaces  
$$C^0 \xrightarrow{d_0} C^1 \xrightarrow{d_1} ... \xrightarrow{d_{n-1}} C^n$$
with $\lbrace c^i \rbrace_{i=0...n}$ and $\lbrace h^i \rbrace_{i=0...n}$ families of bases of the vector spaces $C^i$'s and $H^i$'s, one can define the \textit{torsion} of the based complex $\tor(C^*,\lbrace c^i \rbrace, \lbrace h^i\rbrace)$ as an alternating product of determinants. More precisely, consider the exact sequences
$$0 \to Z^i \to C^i \xrightarrow{d_i} B^{i+1} \to 0$$
$$0 \to B^i \to Z^i \to H^i \to 0$$
that define the vector spaces $B^i$, $Z^i$ and $H^i$.
Pick a system of bases $\lbrace b^i\rbrace$ of the $B^i$'s, first one obtains a basis $b^i \sqcup \bar{h}^i$ of $Z^i$ for any $i$, given by any choice of a section $H^i \to Z^i$. Then any section $B^{i+1} \to C^i$ provides a basis of $C^i$ denoted by $ b^i\sqcup \bar{h}^i \sqcup \bar{b}^{i+1}$. Now compare this new basis with the original basis $c^i$, and take the determinant of the change of basis matrix, denoted by $[b^i\sqcup \bar{h^i} \sqcup \bar{b}^{i+1}: c^i]$. One can show that the alternating product of those determinants does not depend on the lifts, neither on the choice of basis $\lbrace b^i \rbrace$. We define 
 $$\tor(C^*,\lbrace c^i \rbrace, \lbrace h^i \rbrace)= \prod_i [b^i\sqcup \bar{h}^i \sqcup \bar{b}^{i+1}: c^i]^{(-1)^i} \in \mK^*/\lbrace \pm 1 \rbrace$$
 
\begin{remark}
It can be seen as a generalization of the determinant: in particular if the complex is just an isomorphism $(C^0,c^0) \xrightarrow{d_0} (C^1,c^1)$, then the torsion $\tor(C^*,\lbrace c^* \rbrace, \emptyset)$ is nothing but the inverse of the determinant of the invertible matrix of the map $d_0$ in the bases $c^0$ and $c^1$. Note that we defined the torsion up to sign indeterminacy. For the use we will make in this article it makes no difference because we want to study vanishing properties of the torsion, nevertheless we stress out that this sign indeterminacy can be solved in our setting, for instance in \cite{Dub06, DHY09}.
\end{remark}

\subsubsection{The Euler isomorphism}
Given $V$ an $n$-dimensional $\mK$-vector space, its determinant vector space $\det(V)=~\bigwedge^nV$ is defined as its $n$-th exterior power. It is a one dimensional vector space: if $\lbrace v_1, \ldots, v_n\rbrace$ is a basis of $V$, there is an isomorphism $\det(V) \to \mK$ obtained by sending the vector $v_1\wedge \ldots \wedge v_n$ to $1$. In the following, for $L$ a one-dimensional vector space, the notation $L^{\otimes(-1)}$ will denote the dual vector space $\Hom(L,\mK) =L^*$. 
One has the following lemma:
\begin{lemma}\cite[Appendix A, Lemma 5]{GKZ94}
\label{lem:DeterminantIsom}
Let $0 \to A \to B \to C \to 0$ be an exact sequence of vector spaces, then there is a natural isomorphism $$\det(A)\otimes \det(C) \simeq \det(B).$$
\end{lemma}

\begin{definition}
Let $V^*=\bigoplus V^i$ be a finite dimensional graded $\mK$-vector spaces. The \textit{determinant} of $V^*$ is defined by 
$$\det(V^*) = \bigotimes\limits_i \det(V^i)^{\otimes(-1)^i}.$$
\end{definition}
Given a complex $C^*$, the cohomology of this complex is naturally graded, and we have the following proposition that follows easily from Lemma \ref{lem:DeterminantIsom}:
\begin{proposition}\cite[Appendix A, Proposition 3]{GKZ94}
\label{prop:Euler}
There is a natural (Euler) isomorphism $\operatorname{Eu} \colon \det(C^*) \xrightarrow{\sim} \det(H^*(C^*))$.
\end{proposition}

Now fix a based complex $(C^*, c^*)$, where for any $i$, the basis is denoted by $c^i = \lbrace c_1^i, \ldots, c_{n_i}^i \rbrace$. Then we denote by $\bigwedge c^i$ the basis element $c_1^i\wedge \ldots \wedge c_{n_i}^i$ of $\det(C^i)$, and by $\textbf{c} = \bigotimes_i(\bigwedge c^i)^{\otimes(-1)^i}$ the induced basis of the vector space $\det(C^*)$. 
\begin{definition}
\label{defi:Torsion}
The torsion of the based complex $(C^*,c^*)$ is 
$$\tor(C^*,\textbf{c}) = \operatorname{Eu}(\textbf{c}) \in \det(H^*(C^*))/\lbrace\pm1\rbrace.$$
\end{definition}
\begin{remark}
\label{remk:DefinitionsTorsion}
The notation is meaningful: in fact the torsion does not depend on the basis $c^*$ of the complex $C^*$, but only on the basis element $\textbf{c}$ of $\det(C^*)$. Moreover, it coincides with the definition provided in \ref{subsub:Definition} in the following sense:
if $h^*$ is a basis of the graded vector space $H^*(C^*)$, then it defines a basis element $\textbf{h}=\bigotimes_i (\bigwedge h^i)^{\otimes(-1)^i}$, and one can compare $\operatorname{Eu}(\textbf{c})$ with $\textbf{h}$ in $\det(H^*(C^*))$. It provides an element of $\mK$, that we denote by $[\operatorname{Eu}(\textbf{c}):\textbf{h}]$, and we have
$$\tor(C^*,\lbrace c^i\rbrace,\lbrace h^i \rbrace) = [\operatorname{Eu}(\textbf{c}):\textbf{h}].$$
\end{remark}
\subsubsection{The Cayley formula}
\label{subsub:Cayley}
When $C^*$ is an exact complex, a first occurence of a description of the torsion can be found in the seminal work of Cayley in 1848 (see \cite[Appendix B]{GKZ94} where the original text is retranscribed).

Let $(C^*, c^*)$ be a based complex of $\mK$-vector space of the form
$$0 \to (C^0,c^0) \xrightarrow{d_0} \ldots \xrightarrow{d_{r-1}} (C^r,c^r) \to 0.$$ Assume that this complex is exact: it has trivial homology. In particular the one-dimensional vector space $\det(H^*(C^*))$ is canonically isomorphic to $\mK$. We abusively denote by $\tor(C^*, \textbf{c})$ the element of $\mK^*/\lbrace \pm 1\rbrace$ given by $[\tor(C^*, \textbf{c}):1]$, and the equality of Remark \ref{remk:DefinitionsTorsion} $$\tor(C^*, \lbrace c^* \rbrace, \emptyset) = \tor(C^*, \textbf{c})$$ between the two definitions of the torsion holds.

For each index $i$, the basis $c^i = \lbrace c^i_1, \ldots, c^i_{n_i} \rbrace$ can be partitioned into two subsets $c^i_I$ and $c^i_J$ such that $\ker d_i = \langle c^i_I \rangle$. Hence we have $C^i =  \langle c^i_I \rangle \oplus  \langle c^i_J \rangle$ and the map $d_i$ restricts to a linear isomorphism $(d_i)_{I,J}\colon \langle c^i_J \rangle \to \operatorname{im}(d_i) = \langle c^{i+1}_I \rangle$ whose determinant we denote by $\Delta_i$. Of course, those determinants depend on the choices, but it can be shown that their alternating product do not, and we have the proposition:
\begin{proposition}\cite[Appendix A, Theorem 14]{GKZ94}
\label{prop:ProdDet}
With the preceding notations, the following equality holds:
$$\tor(C^*, \textbf{c}) = \prod\limits_{i=0}^{r-1} \Delta_i^{(-1)^{r-i-1}} \in \mK^*/\lbrace \pm1\rbrace.$$
\end{proposition}

\subsection{The Reidemeister torsion form}
\label{sub:TorsionForm}
Recall that $M$ is an irreducible 3-manifold with toral boundary, $Y$ denotes the smooth projective model of a one-dimensional essentially reduced component of $\bar{X}(M)$ of irreducible type, with function field $k(Y)$ and tautological representation $\rho\colon \pi_1(M) \to \SL_2(k(Y))$.
In this subsection we consider the particular case of a the complex $C^*(M, \Adr)$ of twisted cohomology of a 3-manifold $M$. The main goal is to show that in the setting of this article, it defines a rational differential form on the augmented character variety of $M$. 

\medskip
 We want to consider the torsion of the complex $C^*(M, \Adr)$ with coefficients in $\slf_2(k(Y))$. We need to fix bases for the vector spaces $C^i(M, \Adr)$.

\subsubsection{A basis element of $\det(C^*(M, \Adr))$.}
\label{subsub:Basis}
We fix once for all a basis of $\slf_2(k(Y))$ given by $E=~\bsm 0&1\\0&0\esm$, $F=\bsm0&0\\1&0\esm$ and $H=\bsm 1&0\\0&-1\esm$. Now recall that the vector space $C^*(M, \Adr)$ is defined as $C^*(M, \Adr) = \Hom_{\pi_1(M)}(C_*(\widetilde{M}),\slf_2(k(Y)))$: we need to fix a cellular decomposition of $\widetilde{M}$. Since $\partial M \neq \emptyset$, $M$ deformation retracts on a two-dimensional CW-complex: we choose a cell decomposition of this complex $\lbrace e_0^1, \ldots, e_0^{n_0},e_1^1, \ldots, e_1^{n_1}, e_2^1, \ldots, e_2^{n_2} \rbrace$, and we lift it to $\widetilde{M}$. It defines a two-dimensional CW-complex that is a deformation retract of $\widetilde{M}$, we denote the lifts by $\lbrace \widetilde{e}_0^1, \ldots, \widetilde{e}_0^{n_0}, \widetilde{e}_1^1, \ldots, \widetilde{e}_1^{n_1}, \widetilde{e}_2^1, \ldots, \widetilde{e}_2^{n_2} \rbrace$. Now for any $i \in \lbrace 0,1,2\rbrace$, for any $k_i =~1, \ldots, n_i$, for any $\Theta =E,F,H$, one defines $f_i^{k_i,\Theta}$ the element of $C^i(M, \Adr)$ given by $f_i^{k_i,\Theta}(\widetilde{e}_i^{l_j})=\delta_{k_i,l_j}\Theta$. Hence for any $i$, the set $$f_i=\lbrace f_i^{k_i,\Theta}, k_i=1, \ldots, n_i, \Theta=E,F,H \rbrace$$ is a basis of $C^i(M, \Adr).$ We denote by $\textbf{f}$ the basis element of $\det(C^*(M, \Adr))$ given by $\bigotimes_i(\bigwedge f^i)^{\otimes(-1)^i}$.

\subsubsection{Independence of the choices}
Once \textbf{f} is defined, we define the torsion of the complex $C^*(M, \Adr)$ as in Definition \ref{defi:Torsion}: 
$$\tor(C^*(M, \Adr),\textbf{f}) =\operatorname{Eu}(\textbf{f})$$  in the vector space  $\det(H^*(M, \Adr))/\lbrace\pm1\rbrace$. We need to show that it does not depend on the choices.

The main step is the following difficult theorem. It is basically due to Chapman and Cohen, see \cite{Chapman}, \cite{Cohen}.
\begin{theorem}
The torsion does not depend on the choice of the CW-complex for $M$.
\end{theorem}

The remaining part of the invariance is summarized in the following proposition, see for instance \cite[Section 0.2]{Porti97}:
\begin{proposition}
The torsion does not depend on the choice of the basis for $\slf_2(k(Y))$, neither on the choice of the lifts of the cells to $\widetilde{M}$, or on the choice of the tautological representation $\rho\colon \pi_1(M) \to \SL_2(k(Y))$.
\end{proposition}
Consequently, from now on the torsion will denoted by $\tor(M, \Adr)$.

\subsubsection{The torsion form}
Now we must argue why this torsion is in fact an element of $\Omega^1_{k(Y)/k}$. It is enough to prove the following lemma:
\begin{lemma}
\label{lem;TorsionForm}
Let $Y$ be the smooth projective model of a one-dimensional component of irreducible type essentially reduced in the augmented character variety $\bar{X}(M)$, and $\rho\colon \pi_1(M) \to \SL_2(k(Y))$ be a tautological representation.
There is a canonical isomorphism 
$$\det(H^*(M, \Adr)) \simeq \Omega^1_{k(Y)/k}$$
\end{lemma}
\begin{proof}
Recall that the vector space $\det(H^*(M, \Adr))$ is defined as the alternating product $$\det(H^*(M, \Adr)) = \bigotimes_i \det(H^i(M, \Adr))^{\otimes(-1)^i}.$$
We know by Proposition \ref{prop:CohomDimension} that for $i=0$ or $i \ge 2$, the determinant vector spaces $\det(H^i(M, \Adr))$ are naturally isomorphic to $k(Y)$. Since $\dim(Y)=1$, we have that $\dim(H^1(M, \Adr))=1$ hence $\det(H^1(M, \Adr))^*=H^1(M, \Adr)^*$, which is naturally isomorphic to $H_1(M, \Adr)$. We conclude with Corollary \ref{cor:Cotangent} that this is $\Omega^1_{k(Y)/k}$, and the lemma is proved.
\end{proof}
\begin{definition}
\label{defi:TorsionForm}
We define the \textit{Reidemeister torsion form} as $$\tor(M, \Adr) \in \Omega^1_{k(Y)/k}.$$
\end{definition}

\section{The torsion form at finite points}
\label{sec:Finite}
In this section we aim to study the behavior of the torsion form defined in Subsection \ref{sub:TorsionForm} at finite points of the augmented character variety. We will be interested in the vanishing order of this rational differential form. The first striking result is that, although it may \textit{a priori} have poles as well as zeros, the torsion form does not vanish at those finite points. More precisely, we prove the theorem:
\begin{theorem}
\label{theo:Finite}
Let $Y$ be the smooth projective model of a one-dimensional component of irreducible type $\bar{X}$ essentially reduced in the augmented character variety $\bar X(M)$. If $v \in Y$ be a finite point, then the torsion form $\tor(M, \Adr)$ has a pole at $v$ with order equal to the length of the torsion submodule $T_1(M)$ of the $\cO_v$-module $H_1(M, \Adr)_v$. 
\end{theorem}

This torsion submodule will be related with the singularities of the character variety.
There are mainly two possibilities for $v$ to project on a singular character $\chi$ of $X(M)$. 

The first case is that $\chi$ can be a reducible character in $X$, it is necessarily singular since it lies at the intersection of two distinct components: $X$ and the component of reducible type $X(\mZ) \subset X(M)$ that arises from the epimorphism $\pi_1(M) \to \mZ$. We say that it is a singular character of type I. There is a well-know relation (see  \cite{Bur, DR67}) with the Alexander module, in particular such a reducible character $\chi$ induces an eigenvalue $\la_\chi \in k^*$ whose square is a root of the Alexander polynomial.

The second case occurs when $\chi$ is an irreducible singular character (singular of type II). It may be, again, of two different kinds: either it is a proper singularity of the component $X$, either it is an intersection point of $X$ with an other component of irreducible type. To our knowledge, no example of the first situation are known among character varieties of 3-manifolds, while the second situation is illustrated for instance in \cite{Chu17}, see Example \ref{ex:74}.

We can already interpret Theorem \ref{theo:Finite} as follows:
\begin{corollary}
\label{cor:Smooth}
The torsion form 
has no pole or zero at a point $v$ that projects to a smooth character in the character variety $X(M)$.
\end{corollary}
\begin{proof}
Since the length of a module is a non-negative number, the torsion has no zero at finite points $v$ in $Y$. Moreover, it comes from the discussion above that a smooth character $\chi$ in $X$ is in particular irreducible, hence the vector space $H^1(M, \Adbr)$ is isomorphic to the Zariski tangent space $T^\text{Zar}_\chi X(M)$ (\cite{Weil64}). In particular it is a one-dimensional $k$-vector space, so is $H_1(M, \Adbr)$. By the Universal Coefficient Theorem, one has $$H_1(M, \Adbr) \simeq H_1(M, \Adr)_v \otimes_{\cO_v} \cO_v/(t)$$ and it follows from Lemma \ref{lem:Homology} (\ref{3}) and duality that $H_1(M, \Adr)_v \simeq \cO_v$, and $T_1(M) = \lbrace 0 \rbrace$. We conclude that the torsion has no pole at smooth $v$.
\end{proof}

We give the following interpretation of the vanishing order of the torsion at singular points:
\begin{theorem}
\label{theo:Singular}
Assume that $v$ projects to a singular character of type I, such that the associated eigenvalue $\la_\chi$ is a simple root of the Alexander polynomial. Then the torsion form $\tor(M, \Adr)$ has no pole at $v$.

If $v$ projects to a singular character of type II, the order of the pole of the torsion is a (computable) invariant of branch of the singularity, given by the length of the torsion part of the $\cO_v$-module $\Omega_{\bar{B}[M]/k} \otimes_{\bar{B}[M]} \cO_v$. \end{theorem}
\begin{remark}
\label{remk:HigherOrder}
The interested reader can fin a slightly more general statement of the first part of Theorem \ref{theo:Singular} in the author's PhD thesis (\cite[Section 2.3.2]{BenThesis}), that deals with roots of higher order. Note however that there was a mistake in a previous version of this theorem, and that in \cite{BenThesis} one has to invert the  role of poles and zeros at finite points.
Since the result presented there holds under restrictive hypothesis, and that it involves real technical complications, for sake of conciseness  we prefer not to include this generalization to this article, and to focus here on the simplest case of roots of order one. 

\end{remark}

The section is organized as follows: in Subsection \ref{sub:ProofFinite} we prove Theorem \ref{theo:Finite},  in Subsection \ref{sub:Singular I} we prove the first part of Theorem \ref{theo:Singular}, and in Subsection \ref{sub:Singularities} we explain the second part of the statement of Theorem \ref{theo:Singular} and give an example of computation.

\subsection{Proof of Theorem \ref{theo:Finite}}
\label{sub:ProofFinite}
The strategy of the proof is the following: in a first step we produce an acyclic complex, and show that the vanishing order of the torsion of this new complex determines the vanishing order of the torsion form at $v$. Then we compute this vanishing order in terms of the length of $T_1(M)$.

\subsubsection{The cone construction}
We use a well-know construction (see for instance \cite[Appendix A]{GKZ94}) to produce a complex whose torsion will be easier to compute.

We fix a convergent tautological representation $\rho \colon \pi_1(M)\to \SL_2(\cO_v)$, such that its restriction to $\pi_1(\partial M)$ is diagonal. There is a tautological eigenvalue, denoted $\la \colon \pi_1(\partial M) \to~\cO_v^*$ such that for any $\g \in \pi_1(\partial M)$, the matrix $\rho(\g)$ is equal to $\bsm \la(\g)&0\\0&\la^{-1}(\g)\esm$ (see Proposition \ref{prop:TautRepr}).
Recall that we denote by $\rho(\g)_0$ the trace-free matrix $\rho(\g)-\frac 1 2 \Tr(\rho(\g)) I$ in $\slf_2(\cO_v)$.

We define the morphisms of $\cO_v$-modules $\al \colon C^1(M, \Adr)_v \to (\Omega^1_{\cO_v/k})^*$ by
\begin{align*}
&\al(f)\colon \Omega^1_{\cO_v/k} \to \cO_v \\
& d(I_\g \otimes 1) \mapsto \Tr(f([\g])\rho(\g)_0) \text{ for all } \g \in \pi_1(M)\\
& d(1\otimes Z_\g) \mapsto \Tr\left(f([\g])\bsm \la(\g)/2&0\\0&-\la(\g)/2\esm\right) \text{ for all } \g \in \pi_1(\partial M)
\end{align*}
and $\be\colon C^2(M, \Adr)_v \to H^2(M, \Adr)_v \to H^2(\partial M, \Adr)_v \to H^0(\partial M, \Adr)_v^*$ where the first map is the canonical map given by $$0 \to B^2(M, \Adr)_v \to Z^2(M, \Adr)_v \to H^2(M, \Adr)_v \to 0$$ (note that $C^2(M, \Adr)_v=Z^2(M, \Adr)_v$ because $C^3(M, \Adr)_v=0$), the second map is the restriction map induced by the inclusion $\pi_1(\partial M) \to \pi_1(M)$ and the third map is the Poincar\'e duality.

\begin{definition}
Let $R$ be an integral ring with $\mK$ its fraction field. A complex $C^*$ of $R$-modules is \textit{rationally acyclic} if the complex $C^* \otimes_R \mK$ is acyclic. A morphism $\phi\colon C^* \to D^*$ of complexes of $R$-modules is \textit{rationally} a \textit{quasi-isomorphism} if it induces an isomorphism in the cohomology of the rational complexes $C^*\otimes \mK$ and $D^*\otimes \mK$.
\end{definition}

Consider the diagram of complexes 
\begin{center}
\begin{tikzpicture}
\tikzset{node distance=3cm, auto}
\node (A) {$C^0(M,\Ad\circ\rho)_v$};
\node (B) [right of=A] {$C^1(M,\Ad\circ\rho)_v$};
\node (C) [right of=B] {$C^2(M,\Ad\circ\rho)_v$};
\node (D) [above of=A,yshift=-1cm] {$\lbrace 0 \rbrace $};
\node (E) [above of=B,yshift=-1cm] {$(\Omega^1_{\cO_v/k})^*$};
\node (F) [above of=C,yshift=-1cm] {$H^0(\partial M,\Ad\circ\rho)_v^*$};
\draw[->] (A) to node {$d_0$} (B);
\draw[->] (B) to node {$d_1$} (C);
\draw[->] (A) to node {$0$} (D);
\draw[->] (B) to node {$\al$} (E);
\draw[->] (C) to node {$\be$} (F);
\draw[->] (D) to node {$0$} (E);
\draw[->] (E) to node {$0$} (F);
\end{tikzpicture}
\end{center}
We denote by $D^*(M)_v$ the complex with trivial morphisms given by the upper row in this diagram.
We have the following proposition:
\begin{proposition}
\label{prop:Cone}
The morphisms $0$, $\al$ and $\be$  of the diagram above induce a morphism of complexes of $\cO_v$-modules $C^*(M,\Adr)_v \xrightarrow{\phi} D^*(M)_v$. Moreover, the morphism $\phi$ is rationally a quasi-isomorphism.
\end{proposition}
\begin{proof}
To prove the first claim we need to show that this diagram commutes. It is clear from the definition of the morphism $\be$ that the composition $\be \circ d_1$ is zero.
Let $\zeta$ be a $0$-cochain in $C^0(M, \Ad\circ\rho)_v$. For any $\g$ in $\pi_1(M)$, we have 
$$\al(d_0\zeta) \left(d(I_\g \otimes 1)\right) = \Tr(d_0\zeta(\g)\rho_0(\g)).$$
But $d_0\zeta(\g) = \rho(\g)\zeta \rho^{-1}(\g) - \zeta$ and for any $\g$ in $\pi_1(M)$, $\rho(\g)\rho_0(\g) = \rho_0(\g)\rho(\g)$, hence 
\begin{align*}
\Tr(d_0\zeta(\g)\rho_0(\g)) = \Tr(\rho(\g)\zeta\rho^{-1}(\g)\rho_0(\g)) - \Tr(\zeta\rho_0(\g))\\
=\Tr(\rho(\g)\zeta\rho_0(\g)\rho^{-1}(\g))-\Tr(\zeta\rho_0(\g))=0.
\end{align*}
For $\gamma \in \pi_1(\partial M)$, a similar computation shows that $\al(d_0\zeta) \left(d(1\otimes dZ_\gamma)\right)=0$ and it proves the first statement of the proposition.

Now we prove that the two complexes $C^*(M, \Adr)$ and $D^*(M)_v \otimes~k(Y)=~D^*(M)$ are quasi-isomorphic.
We start by noticing that the two vector spaces $H^0(M, \Adr)$ and $H^0(D^*(M))$ are trivial. Then it follows from \cite[Chapter 6, Proposition 1.8, (c)]{QL06} that $(\Omega^1_{\cO_v/k})^* \otimes k(Y) \simeq (\Omega^1_{k(Y)/k})^*$, hence the map induced by $\al$ in cohomology $H^1(M, \Adr) \xrightarrow{\al_*} (\Omega^1_{k(Y)/k})^*$ is nothing but the dual map of the isomorphism constructed in the proof of Corollary \ref{cor:Cotangent}. Finally $\be$ induces the isomorphism $H^2(M, \Adr) \simeq H^0(\partial M, \Adr)^*$ (see the proof of Proposition \ref{prop:CohomDimension}), and the proposition is proved.
\end{proof}

Now we use a standard construction, called \textit{cone} of the morphism $\phi$ (see for instance \cite[Appendix A]{GKZ94}). We define $D^{*-1}_v$ as the complex $D^*_v$ but with the numbering shifted by one, namely the $0$-th $\cO_v$-module is the trivial module (added) with the trivial boundary map, then the first $\cO_v$-module is the zero-th $\cO_v$-module from $D^*_v$, and so on. The complex $\operatorname{Cone}(\phi)$ is the complex of $\cO_v$-modules $D^{*-1}(M)_v \oplus C^{*}(M, \Adr)_v$ given by 
$$C^0(M, \Adr)_v \xrightarrow{d_0} C^1(M, \Adr)_v \xrightarrow{d_1 \oplus \al} C^2(M, \Adr)_v \oplus (\Omega^1_{\cO_v/k})^*\xrightarrow{\beta \oplus 0} H^0(\partial M, \Adr)_v^* $$
The following lemma follows immediately from Proposition \ref{prop:Cone}:
\begin{lemma}
\label{lem:ConeAcyclic}
The complex of $\cO_v$-modules $\operatorname{Cone}(\phi)$ is rationally exact.
\end{lemma}
We want to compare the torsion of the acyclic rational complex $\operatorname{Cone}(\phi) \otimes k(Y)$ with the torsion of the complex $C^*(M, \Adr)$. To do so, the first step is to fix bases of the complex $\operatorname{Cone}(\phi) \otimes k(Y)$. Since we already fixed bases of $C^*(M, \Adr)$, we only need a basis of $(\Omega^1_{k(Y)/k})^*$ and of $H^0(\partial M, \Adr)^*$. Recall that $t$ is an element of $k(Y)$ with valuation $v(t)$ equal to $1$. In particular $dt$ is not zero in $\Omega^1_{k(Y)/k}$, and we take it as a basis element. Since the tautological representation $\rho_{\partial M}$ is diagonal when restricted to $\pi_1(\partial M)$, we take the matrix $H = \bsm 1&0\\0&-1\esm$ as a basis of $H^0(\partial M, \Adr)$, and we denote by $H^*$ the dual basis of  $H^0(\partial M, \Adr)^*$.
\begin{remark}
\label{remk:Vanishing}
This choice is consistent with the identification $\det(H^2(M, \Adr)) \simeq k(Y)$ of Propostion \ref{prop:CohomDimension}. In particular with those choices of bases, the torsion $\tor(M, \Adr)$ can be written as $P\cdot dt \otimes H^*$, for some $P\in k(Y)$. What we need to compute is the vanishing order of the function $P$, and to do so, we relate it with the torsion of the complex $\operatorname{Cone}(\phi) \otimes k(Y)$.
\end{remark}
We denote by $\mathcal{H}$ the long exact sequence in cohomology induced by the exact sequence of $k(Y)$-vector spaces $$0 \to D^{*-1}(M) \to \operatorname{Cone}(\phi) \otimes k(Y) \to C^{*}(M, \Adr) \to 0.$$ Since the complex $\operatorname{Cone}(\phi) \otimes k(Y)$ is acyclic it splits into two isomorphisms
\begin{equation}
\label{equa:LongSequence}
H^1(M, \Adr) \xrightarrow{\al_*} (\Omega^1_{k(Y)/k})^* \text{ and }\ H^2(M, \Adr) \xrightarrow{\be_*} H^0(\partial M, \Adr)^*
\end{equation}

To relate the torsion $\tor(M, \Adr)$ with the torsion of the rational complex $\operatorname{Cone}(\phi) \otimes k(Y)$, we use the well-known multiplicativity of the torsion, see \cite[Theorem 3.2]{Mil66}.
\begin{proposition}
\label{prop:Multiplicativity}
Let $0 \to (A^*,\textbf{a}) \to (B^*,\textbf{b}) \to (C^*, \textbf{c}) \to 0$ be an exact sequence of based complexes of vector spaces, with bases of their homology $\textbf{h}_a, \textbf{h}_b$ and $\textbf{h}_c$. Let $\mathcal{H}$ be the induced long exact sequence in cohomology, seen as an acyclic based complex. Then the following equality holds:
\begin{equation}
\label{equa:Mult}
\tor(B^*, \textbf{b},\textbf{h}_b,)=\tor(A^*, \textbf{a},\textbf{h}_a,) \tor(C^*, \textbf{c},\textbf{h}_c,) \tor(\mathcal{H},\textbf{h}_a,\textbf{h}_b,\textbf{h}_c)
\end{equation}
\end{proposition}

In particular equation (\ref{equa:Mult}) turns into 
\begin{equation}
\label{equa:TorsionCone}
\tor(\operatorname{Cone}(\phi) \otimes k(Y), \textbf{f} \sqcup dt \otimes H^*, \emptyset) = \frac{\tor(C^*(M, \Adr), \textbf{f},dt \otimes H^*)}{\tor(D^*(M),dt \otimes H^*,dt \otimes H^*)} \tor(\mathcal{H},dt \otimes H^*, \emptyset)
\end{equation}
where $\textbf{f}$ is the basis element in $\det C^*(M, \Adr)$ fixed in Section \ref{sub:TorsionForm}.

Since $\tor(D^*(M),dt \otimes H^*,dt \otimes H^*) = 1$ and $\tor(C^*(M, \Adr), \textbf{f},dt \otimes H^*) = P$, equation (\ref{equa:TorsionCone}) becomes
\begin{equation}
\label{equa:TorsionCone2}
\tor(\operatorname{Cone}(\phi) \otimes k(Y), \textbf{f} \sqcup dt \otimes H^*, \emptyset) =P \, {\tor(\mathcal{H},dt \otimes H^*)}
\end{equation}
and we are led to compute the numerator of the right-hand side term in (\ref{equa:TorsionCone2}).
\subsubsection{The torsion of the long exact sequence}
We prove the following lemma:
\begin{lemma}
\label{lem:TorsionLong}
The torsion of the long exact sequence $\mathcal{H}$ is an invertible element of $\cO_v$. In particular the vanishing order of $\tor(\operatorname{Cone}(\phi) \otimes k(Y), \textbf{f} \sqcup dt \otimes H^*, \emptyset)$ at $v$ is given by $v(P)$.
\end{lemma}
   
To see this, observe that computing the torsion of this sequence turns out to compute the valuation of the (inverses of) the determinants of the isomorphisms $\al_*$ and $\be_*$ of \ref{equa:LongSequence}. We shortly explain how to do such a computation:
those isomorphisms descend to $\cO_v$-modules homomorphisms $\al_{*,v} \colon H^1(M, \Adr)_v \to (\Omega^1_{\cO_v/k})^*$ and $\be_{*,v}\colon H^2(M, \Adr)_v \to H^0(\partial M, \Adr)_v^*$. By Lemma \ref{lem:Homology}, $H^1(M, \Adr)_v$ is a free module of rank one while $H^2(M, \Adr)_v \simeq \cO_v \oplus T_1(M)$, and we may rather consider the homomorphisms of free $\cO_v$-modules $\al_{*,v}$ and $\be_{*,v}\colon H^2(M, \Adr)_v/T_1(M) \to H^0(\partial M, \Adr)_v^*$. Since those are rationally isomorphisms, it is clear that they are one-one, and it is follows from the definition that the valuation $v(\det(\al_*))$ (respectively $\be_*$) is nothing but the length of the cokernel $\coker(\al_{*,v})$ (resp. $\coker(\be_{*,v})$).
It can be seen as a corollary of the following theorem which will be very useful along this paper:
\begin{theorem}\cite[Appendix A, Theorem 30]{GKZ94},\cite[Theorem 4.7]{Tur01}
\label{theo:Russes}
Let $(C^*, \textbf{c})$ be a rationally exact based complex of free $\cO_v$-modules. Then the valuation of the torsion of the rational complex $C^* \otimes k(Y)$ can be computed as follows:
$$v(\tor(C^* \otimes k(Y),\textbf{c}, \emptyset))= \sum_k (-1)^k \operatorname{length}(H^k(C^*)).$$
\end{theorem}
It follows from the discussion above that Lemma \ref{lem:TorsionLong} will be proven if we prove that both $\al_{*,v}$ and $\be_{*,v}$ are surjective homomorphisms. It is the aim of the two upcoming lemmas.
\begin{lemma}
\label{lem:Onto1}
The homomorphism $\al_{*,v} \colon H^1(M, \Adr)_v \to (\Omega^1_{\cO_v/k})^*$ is onto.
\end{lemma}
\begin{proof}
We construct a right-section $s \colon (\Omega^1_{\cO_v/k})^* \to H^1(M, \Adr)_v$ for $\alpha_{*,v}$ as follows. Let us fix a linear form $\theta\colon \Omega^1_{\cO_v/k} \to \cO_v$. 

By the universal property of the module of differential forms $\Omega^1_{\cO_v/k}$ (\cite[Definition 1.2, Chapter 6]{QL06}), such a linear form corresponds uniquely to a $k$-derivation still denoted by $\theta \colon \cO_v \to \cO_v$. Given a matrix in $\SL_2(\cO_v)$, one can apply the derivation $\theta$ to each of its entry, it yields a map (still denoted by)  $\theta\colon \SL_2(\cO_v) \to \SL_2(\cO_v)$. Explicitly, this map is given by $\theta(d(I_\gamma \otimes 1)) = \Tr \theta(\rho(\gamma))$ for $\gamma \in \pi_1(M)$, and by $\theta (d(1\otimes Z_\gamma))= \Tr \theta (\lambda(\gamma))$ for $\gamma \in \pi_1(\partial M)$.

Now we define $s(\theta)=s_\theta$ in $H^1(M, \Adr)_v$ by the formula
$$s_\theta \colon \gamma \mapsto \theta(\rho(\g)) \rho(\g)^{-1}$$
and we need to check that it defines an element of $Z^1(M, \Adr)_v$, namely that for any $\gamma$, the matrix $s_\theta(\gamma)$ is trace-free, and that it satisfies the cocyle relation.
\begin{claim}
The matrix $s_\theta(\gamma)$ is trace-free for any $\gamma$ and satisfies the cocycle relation $s(\theta)(\g\de) = s(\theta)(\g) + \Adr(\g) (s(\theta)(\de))$ for any $\gamma, \delta$.
\end{claim}
\begin{proof}[Proof of the claim]
Using the Leibniz rule, a direct computation shows that for any two by two matrix $A$, one has
\begin{equation}
\label{equa:Section}
\Tr(\theta(A)A^{-1}) = \theta(\det A)
\end{equation}
and it yields $\Tr(s_\theta(\gamma)) = \theta(\det \rho(\g))=0$, the last equality being clear since $\theta(1)=0$ by definition. The cocycle relation can be also checked by direct computation.
\end{proof}

To conclude the proof of Lemma \ref{lem:Onto1}, we need to see that the composition $\al_{*,v} \circ s \colon (\Omega^1_{\cO_v/k})^* \to (\Omega^1_{\cO_v/k})^*$ is the identity map. Pick $\theta \in (\Omega^1_{\cO_v/k})^*$, then for any $\gamma$ in $\pi_1(M)$ we have $$\al_{*,v}(s_\theta(\gamma)) =\Tr(s(\theta)(\g)\rho_0(\g)) = \Tr(\theta(\rho(\g)) \rho(\g)^{-1}\rho_0(\g)) = \Tr(\theta(\rho(\g))),$$ the last equality being consequence from (\ref{equa:Section}). Considering again $\theta$ as a linear form, the latter is nothing but $\theta(d(I_\g \otimes 1))$ and it completes the proof of the lemma.
\end{proof}

\begin{lemma}
\label{lem:Onto2}
The homomorphism $\be_{*,v} \colon H^2(M, \Adr)_v \to H^0(\partial M, \Adr)_v^*$ is onto.
\end{lemma}
\begin{proof}
Recall that $\be_{*,v}$ is by construction the composition of the morphism induced by inclusion $H^2(M, \Adr)_v \to H^2(\partial M, \Adr)_v$ with the Poincar\'e duality $H^2(M, \Adr)_v \simeq H^0(M, \Adr)_v^*$. Using the long exact sequence in cohomology of the pair $(M, \partial M)$, one sees that the first morphism is onto since the $\cO_v$-module $H^3(M, \partial M, \Adr)_v \simeq H_0(M, \Adr)_v$ is trivial by Lemma \ref{lem:H0}, and it is clear for the second one.
\end{proof}
\begin{proof}[Proof of Lemma \ref{lem:TorsionLong}]
It follows from the discussion after Lemma \ref{lem:TorsionLong} and from Lemmas \ref{lem:Onto1} and \ref{lem:Onto2} that the determinants of the morphisms of $k(Y)$-vector spaces $\al_*$ and $\be_*$ of (\ref{equa:LongSequence}) are invertible elements of $\cO_v$ and the lemma follows.
\end{proof}

\subsubsection{Proof of Theorem \ref{theo:Finite}}
Now we are ready to prove Theorem \ref{theo:Finite}: the vanishing order of $\tor(M, \Adr))$ is given by $v(P)$ (Remark \ref{remk:Vanishing}) which turns out to be equal to $v(\tor((\operatorname{Cone}(\phi) \otimes k(Y), \textbf{f} \sqcup dt \otimes H^*, \emptyset)))$ by Lemma \ref{lem:TorsionLong}. Let us compute the term $v(\tor((\operatorname{Cone}(\phi) \otimes k(Y), \textbf{f} \sqcup dt \otimes H^*, \emptyset)))$. Since the complex $\operatorname{Cone}(\phi)$ is a rationally acyclic complex of free $\cO_v$-modules, we shall use Theorem \ref{theo:Russes}. We obtain the following equality:
$$v(\tor((\operatorname{Cone}(\phi) \otimes k(Y), \textbf{f} \sqcup dt \otimes H^*, \emptyset)))=\sum_k (-1)^k \operatorname{length}(H^k(\operatorname{Cone}(\phi)))$$
To compute the cohomology of the complex $\operatorname{Cone}(\phi)$, we use the long exact sequence in cohomology $\mathcal{H}(\operatorname{Cone}(\phi))$ induced by the exact sequence of $\cO_v$-modules $$0 \to D^*(M)_v \to \operatorname{Cone}(\phi) \to C^{*+1}(M, \Adr)_v \to 0,$$ we obtain
\begin{align}
\label{equa:LongExactSeq}
0 \to &H^0(M, \Adr)_v \to H^0(D^*(M)) \to H^0(\operatorname{Cone}(\phi)) \to H^1(M, \Adr)_v \to H^1(D^*(M)) \\
& \to H^1(\operatorname{Cone}(\phi)) \to H^2(M, \Adr)_v \to H^2(D^*(M)) \to H^2(\operatorname{Cone}(\phi)) \to 0 \nonumber
\end{align}
But $H^0(M, \Adr)_v = \lbrace 0 \rbrace$ by Lemma \ref{lem:H0}, $H^0(D^*(M))$ is obviously trivial too, $H^1(D^*(M))=(\Omega^1_{\cO_v/k})^*$ and $H^2(D^*(M))= H^0(\partial M, \Adr)_v^*$, hence (\ref{equa:LongExactSeq}) becomes
\begin{align*}
\label{equa:LongExact2}
0 \to H^0(\operatorname{Cone}(\phi)) &\to H^1(M, \Adr)_v \xrightarrow{\al_{*,v}} (\Omega^1_{\cO_v/k})^* \to H^1(\operatorname{Cone}(\phi))\\
& \to H^2(M, \Adr)_v \xrightarrow{\be_{*,v}} H^0(\partial M, \Adr)_v^* \to H^2(\operatorname{Cone}(\phi)) \to 0
\end{align*}
By Lemma \ref{lem:Homology} (\ref{3}) combined with Lemma \ref{lem:Onto1}, it comes that $\al_{*,v}$ is an isomorphism and $H^0(\operatorname{Cone}(\phi))=\lbrace 0 \rbrace$. Then Lemma \ref{lem:Homology} (\ref{5}) with Lemma \ref{lem:Onto2} imply that $H^1(\operatorname{Cone}(\phi)) \simeq T_1(M)$ and $H^2(\operatorname{Cone}(\phi))=\lbrace 0 \rbrace$.

We conclude that $v(\tor((\operatorname{Cone}(\phi) \otimes k(Y), \textbf{f} \sqcup dt \otimes H^*, \emptyset)))= -\operatorname{length}(T_1(M))$, hence we have $v(P) =- \operatorname{length}(T_1(M))$ and the theorem follows from the equality $\tor(M, \Adr) = P\cdot dt$.

\subsection{Reducible characters and Alexander polynomial}
\label{sub:Singular I}
In this subsection we assume that $v$ projects on a reducible character $\chi$ in $X$ (a singular character of type I), and we show the first part of Theorem \ref{theo:Singular}, which we recall as the following proposition:
\begin{proposition}[First part of Theorem \ref{theo:Singular}]
\label{prop:ReducChar}
Assume that a finite point $v$ in $Y$ projects on a reducible character $\chi$ in the character variety $X(M)$, with eigenvalue $\lambda_\chi \in k^*$. If $\lambda_\chi$ is a simple root of the Alexander polynomial $\Delta_M$, then the torsion form has no pole at $\chi$.
\end{proposition}

Recall that we assume that the first Betti number of $M$ is equal to 1, hence there is a unique (up to inversion) abelianization epimorphism $\varphi\colon \pi_1(M) \to \mZ$. 
Any reducible representation $\al\colon \pi_1(M) \to \SL_2(k)$ can be written, for any $\g$ in $\pi_1(M)$, as
\begin{equation}
\label{equa:AlphaLambda}
\al_\la(\g) = \bsm \la^{\varphi(\g)} & * \\0 & \la^{-\varphi(\g)}\esm
\end{equation}
where $\la$ is an element of $k^*$.
\begin{remark}
Recall that in Proposition \ref{prop:TautRepr} we have defined a tautological eigenvalue $\la_Y \colon \pi_1(\partial M)\to k(Y)$. It turns out that for a finite $v \in Y$, for any $\g$ in $\pi_1(\partial M)$, the function $\la_Y(\g)$ lies in $\cO_v$. In this case number $\la^{\varphi(\g)}$ is nothing but the residual eigenvalue $\bar{\la}_Y(\g) = \ev(\la_Y(\g)) \in k^*$.
\end{remark}
Since the eigenvalue $\la$ does not depend on the choice of the representation $\al_\la$ but only on the character $\chi$, we say that it is \textit{the eigenvalue associated to} $\chi$.

\medskip

The following theorem has a long story, it is originally due to Burde (\cite{Bur}) and de Rham (\cite{DR67}), see also \cite{HPS01}, \cite[Theorem 1.4.6]{BenThesis}, \cite[Theorem 2.4]{BC18} for a more recent treatment and generalizations.
\begin{theorem}
Assume that $\lambda \neq \pm 1$, then there exists a reducible, non-abelian representation $\al_\la$ of the form (\ref{equa:AlphaLambda}) if and only if $\la^2$ is a root of the Alexander Polynomial.
\end{theorem}
By Lemma \ref{lem:NonCentral} no reducible character in a component of irreducible type is a central character, hene we obtain the following statement:
\begin{corollary}
Let $\chi$ be a reducible character in a component of irreducible type $X$ of the character variety $X(M)$. Then the square $\la_\chi^2$ of the associated eigenvalue is a root of the Alexander polynomial.
\end{corollary}

\subsubsection{Proof of Proposition \ref{prop:ReducChar}}
To prove Proposiion \ref{prop:ReducChar}, one may consider the $\cO_v$-module $H_1(M, \Adr)$ and show that it is torsion free under the hypothesis that $\lambda^2_\chi$ is a simple root of the Alexander polynomial.

It wiil follow promptly from the lemma:
\begin{lemma}
\label{lem:OneDimRed}
The residual $k$-vector space $H^1(M, \Adbr)$ is one-dimensional.
\end{lemma}

First we deduce the proof of proposition \ref{prop:ReducChar} from Lemma \ref{lem:OneDimRed}, and then we will prove this lemma in the next subsection.
\begin{proof}[Proof of Proposition \ref{prop:ReducChar}]
As already mentioned in this article, the rational $k(Y)$-vector space $H_1(M, \Adr) \simeq H_1(M, \Adr)_v \otimes_{\cO_v} k(Y)$ is one-dimensional, since it is isomorphic to the space of rational one-forms on the one-dimensional variety $Y$. So the $\cO_v$-module $H_1(M, \Adr)_v$ has the form
\begin{equation*}
 \label{equa:Numberl}
H_1(M, \Adr)_v=\cO_v \oplus T_1(M)
\end{equation*}
where $T_1(M)$ is the torsion part of $H_1(M, \Adr)_v$.
Now using the Universal Coefficients Theorem, one gets that
$$H_1(M, \Adr)_v \otimes_{\cO_v} \cO_v/(t) \simeq H_1(M, \Adbr) \simeq H^1(M, \Adbr)$$
is  one-dimensional as a $k$-vector space, hence $T_1(M) = \lbrace 0 \rbrace$.
The proposition follows from Theorem~\ref{theo:Finite}.
\end{proof}
\subsubsection{Proof of Lemma \ref{lem:OneDimRed}}
\label{subsub:step1}
Now we prove the technical lemma \ref{lem:OneDimRed}. The strategy is to make a "d\'evissage" of the vector space $H^1(M, \Adbr)$, and ultimately to use a cohomological interpretation of the fact that $\lambda^2_\chi$ is a simple root of the Alexander polynomial.

We denote by $\varphi_\la\colon \pi_1(M) \to k^*$ the homomorphism $\varphi_\la(\g) = \la^{\varphi(\g)}$, where $\varphi\colon \pi_1(M) \to \mZ$ is the abelianization homomorphism. It makes the group $\pi_1(M)$ acts
 on $k$ by multiplication and we write $\varphi_\lambda$ to emphasize the $\mZ[\pi_1(M)]$-module structure on the field $k$. We denote by $C^*(M, \varphi_\la)$ the twisted complex with coefficients $\varphi_\lambda$, and $H^*(M, \varphi_\la)$ the corresponding cohomology. Note that $H^0(M, \varphi_\la)$ is trivial for any $\la \neq 1$.

By Lemmas \ref{lem:NonCentral} and \ref{lem:NotAbelian} we can fix a convergent tautological representation $\rho\colon \pi_1(M) \to \SL_2(\cO_v)$ that is not residually abelian, so that the residual representation $\brho\colon \pi_1(M)\to \SL_2(k)$ has the form~(\ref{equa:AlphaLambda}). 

More precisely, we have the following well-known lemma that describes the representation $\brho$, see for instance \cite[Section 1.4.2]{BenThesis}, or \cite[Lemma 2.17]{BC18} for a proof.
\begin{lemma}[\cite{DR67}]
\label{lem:Cocycleu}
There is a non-trivial cocycle $u \in H^1(M,  \varphi_{\la^{2}})$ such that the residual representation is of the form $$\brho(\g) = \bm \varphi_\la(\g) & \varphi_\la(\g)^{-1} u(\g)\\0& \varphi_\la(\g)^{-1}\ema$$
for any $\g$ in $\pi_1(M)$
\end{lemma}
We need to fix some notations. The 3-dimensional representation $\Adbr$ can be computed in the basis $\lbrace E,H,F \rbrace$, it has the following form:
\begin{equation}
\label{equa:Adjoint}
\Adbr(\g) = 
\bm \varphi_{\la^2}(\g)& -2u(\g)& -\varphi_{\la^{-2}}(\g)u^2(\g)\\
0 & 1 & \varphi_{\la^{-2}}(\g)u(\g)\\
0 & 0 & \varphi_{\la^{-2}}(\g)\ema
\end{equation}
for any $\g$ in $\pi_1(M)$.
We denote by $\Adbr_{3,3}$ the sub-representation of $\Adbr$ obtained by deleting the third row and the third column of the matrix $\Adbr(\g)$. It acts on $k^2_{\Adbr_{3,3}}= k.E \oplus k.H$. In other words, for any $\gamma$ in $\pi_1(M)$ we have
\begin{equation}
\label{equa:Adjoint33}
\Adbr_{3,3}(\g) = 
\bm \varphi_{\la^2}(\g)& -2u(\g)\\
0 & 1 \ema
\end{equation}
Hence the $\pi_1(M)$-module $\slf_2(k)_{\Adbr}$ splits into the following exact sequence of $\pi_1(M)$-modules, where the induced action is indicated as a subscript:
\begin{equation}
\label{equa:ModuleSplit}
0 \to k^2_{\Adbr_{3,3}} \to \slf_2(k)_{\Adbr} \to \varphi_{\lambda^{-2}} \to 0
\end{equation}
\begin{lemma}
The sequence (\ref{equa:ModuleSplit}) induces the long exact sequence of $k$-vector spaces in cohomology:
\begin{align}
\label{equa:ModulesLong}
0 \to & H^1(M, \Adbr_{3,3}) \to H^1(M, \Adbr) \to H^1(M, \varphi_{\la^{-2}}) \to \\
&H^2(M, \Adbr_{3,3}) \to H^2(M, \Adbr)\to H^2(M, \varphi_{\la^{-2}}) \to 0 \nonumber
\end{align}
\end{lemma}
\begin{proof}
Since $\brho$ is not abelian $\lambda\neq \pm1$ so that $H^0(M, \varphi_{\lambda^{-2}})$ is trivial and the lemma follows.
\end{proof}

Now we are led to study the cohomological complex $C^*(M, \Adbr_{3,3})$. Using (\ref{equa:Adjoint33}), one sees that the $\pi_1(M)$-module $k^2_{\Adbr_{3,3}}$ splits as
\begin{equation*}
0 \to \varphi_{\lambda^2} \to k^2_{\Adbr_{3,3}} \to k\to 0
\end{equation*}
This splitting induces the long exact sequence
\begin{align}
\label{equa:LongSplit2}
0 \to H^0(M) \to H^1(M,  \varphi_{\la^2}) \to H^1(M, \Adbr_{3,3}) \\ \nonumber \to H^1(M) \xrightarrow{\de} H^2(M,  \varphi_{\la^2}) \to H^2(M, \Adbr_{3,3}) \to 0
\end{align}
where $H^i(M)$ is the usual cohomology with non-twisted coefficients $k$. In particular the vector spaces $H^0(M)$ and $H^1(M)$ are one-dimensional since $b_1(M)=1$, and $H^2(M)$ is trivial. 

Now we use the hypothesis that we consider a simple root $\lambda^2_\chi$.
\begin{lemma}
\label{lem:H1Alex}
If $\lambda^2$ is a simple root of the Alexander polynomial, then the $k$-vector space $H^1(M, \varphi_{\lambda^2})$ is one-dimensional.
\end{lemma}
\begin{proof}
This follows directly from well-known facts from Alexander modules theory. In a nutshell, since the Alexander polynomial is symmetric, it can be written $\Delta_M(t) = (t-\lambda^{-2}) P(t)$, and $P(\lambda^{-2}) \neq 0$. Now it implies that the Alexander module has $(t-\lambda^{-2})$-torsion exactly of the form $\cO_v/(t-\lambda^{-2})$, and the lemma follows from the Universal Coefficients Theorem. See for instance \cite[Proof of Lemma 2.8]{BC18} for more details.
\end{proof}
\begin{lemma}
\label{lem:Acyclic33}
If $\lambda^2$ is a simple root of the Alexander polynomial, then the complex $C^*(M, \Adbr_{3,3})$ is acyclic.
\end{lemma}
\begin{proof}
Using Lemma \ref{lem:H1Alex}, we deduce that the first homomorphism in the sequence (\ref{equa:LongSplit2}) $H^0(M) \to H^1(M, \varphi_{\lambda^2})$ is an isomorphism. Hence we are led to consider the homomorphism $\delta \colon H^1(M) \to H^2(M, \varphi_{\lambda^2})$. By a simple diagram chasing, one can compute it explicitely: 
\begin{align*}
\delta \colon H^1(M) &\to H^2(M, \varphi_{\lambda^2})\\
\phi &\mapsto \phi \cup u
\end{align*}
where $u \in H^1(M, \varphi_{\lambda^2})$ is the cocycle given by Lemma \ref{lem:Cocycleu}, and $\cdot \cup \cdot \colon H^1(M) \times H^1(M, \varphi_{\lambda^2}) \to H^2(M, \varphi_{\lambda^2})$ is the cup-product (\cite[Chapter V]{Bro}).
Now it is proved in \cite[Corollary 6.7]{HP15} that this cup product is non-trivial, in particular the map $\delta$ in injective, so that $H^1(M, \Adbr_{3,3})$ vanishes, and so does $H^2(M, \Adbr_{3,3})$ since $\chi(M)=0$.
\end{proof}
Now we are ready for the proof of Lemma \ref{lem:OneDimRed}.
\begin{proof}[Proof of Lemma \ref{lem:OneDimRed}]
Inserting the result of Lemma \ref{lem:Acyclic33} in the sequence (\ref{equa:ModulesLong}), one gets the isomorphism
$$H^1(M, \Adbr) \simeq H^1(M, \varphi_{\lambda^{-2}})$$
and the latter is one-dimensional by Lemma \ref{lem:H1Alex} and Blanchfield duality.
\end{proof}

\subsection{Singularities of algebraic curves}
\label{sub:Singularities}
In this subsection we explain the second part of Theorem \ref{theo:Singular}.

\medskip

Recall that we picked $\bar{X}$ a one-dimensional component of the augmented character variety $\bar{X}(M)$, and $Y$ its smooth projective model: it comes with a birational morphism $\nu \colon Y \to \bar{X}$. More precisely, $\nu$ is defined at any finite point $v$ of $Y$, and it is an isomorphism in restriction to the smooth locus of $\bar{X}$. However, $\bar{X}$ might be singular.

Whatever $\nu(v)$ is singular or smooth, if it projects to an irreducible character $\chi$ of $X$, then we have (Lemma \ref{lem:Homology}, (\ref{2}))
\begin{equation*}
H_1(M, \Adr)_v \simeq \Omega^1_{\bar{B}[M]/k} \otimes \cO_v
\end{equation*}
and the torsion form $\tor(M, \Adr)$ has a pole as prescribed by the length of the torsion part of the $\cO_v$-module $\Omega^1_{\bar{B}[M]/k} \otimes \cO_v$.

One may wonder about the interpretation of this $\cO_v$-module. As we already proved, this module is torsion-free if $\nu(v)$ is a smooth point of $\bar{X}$. On the other hand, to our knowledge it is not known if the converse statement is true (see \cite{Ber94} for a survey on a more general question).

A singular point $x$ in the curve $\bar{X}$ may have several pre-images in $Y$, each of them corresponds to a \textit{branch} around $x$ in $X$.
For instance the cusp $(0,0) \in \lbrace X^2-Y^3=0 \rbrace$ has only one branch, the normal double point $(0,0) \in \lbrace XY=0 \rbrace$ has two branches...

In the following, since the module $\Omega^1_{\bar{B}[M]/k} \otimes \cO_v$ carries only local informations, we study its behavior for the family of plane singularities $(0,0) \in \lbrace X^p-Y^q=0 \rbrace, p<q$. We denote by $C$ the curve $\lbrace (x,y) \in k^2 \vert \ x^p-y^q = 0 \rbrace,$ and by $x$ the point $(0,0)$ in $C$. Let $n= \operatorname{gcd}(p,q)$, $p'=\frac p n$, $q'=\frac q n$. The morphism of $k$-algebras 
\begin{align*}
\nu \colon k[X,Y]/(X^p-Y^q) &\to k[S]\\
X \mapsto S^{q'}, Y\mapsto S^{p'}
\end{align*}
induces a birational morphism $\nu\colon \mathbb{A}^1_k \to C$. The preimage of $x$ is the point $0$ in $\mathbb{A}^1_k$. The module of K\"ahler differentials at $x$ is $\Omega_{\cO_x/k} \simeq \cO_x dX \oplus \cO_x dY /(pX^{p-1}dX-qY^{q-1}dY),$ hence its pull-back by $\nu$ is $\Omega_{\cO_x/k} \otimes_{\cO_x} \cO_0 = \cO_0 dX \oplus \cO_0 dY/(pS^{q'(p-1)}dX-qS^{p'(q-1)}dY)$. Consider the morphism induced by $\nu$:
\begin{align*}
\Omega_{\cO_x/k} &\otimes \cO_0 \to \Omega_{\cO_0/k}\\
dX &\mapsto q'S^{q'-1}dS\\
dY &\mapsto p'S^{p'-1}dS
\end{align*}
Since $0$ is smooth, the right-hand side $\Omega_{\cO_0}$ is a free $\cO_0$-module of rank one, hence the kernel of this morphism is the torsion part of $\Omega_{\cO_x/k} \otimes \cO_0$. It is generated by $(p'dX -q'S^{q'-p'}dY)$ and its annihilator is $(nS^{q'(p-1)})$. In particular the length of the torsion part of the module $\Omega_{\cO_x/k} \otimes \cO_0$ is $q'(p-1)$.

\section{Examples and computations}
\label{sec:Examples}
In this section we use the computations of \cite{Porti97,Dub06,DHY09} to give explicit examples of the torsion form. In Subsection \ref{sub:Formula} we relate the torsion form with the work of the sub-mentioned authors. In Subsection \ref{sub:Compute} we use this relation to give polynomial expressions for the torsions of four simple knots exteriors in $\mS^3$.

\subsection{A comparison formula}
\label{sub:Formula}
In this subsection we give a formula that relates the torsion function defined in \cite{Porti97} with the differential form of this article. 

\medskip

Let $\mu$ be a curve in $\pi_1(\partial M)$, we denote by $\mT_\mu$ the torsion function defined in \cite{Porti97}. 
We prove the following proposition:
\begin{proposition}
\label{prop:Formula}
Let $X$ be a component of irreducible type of the character variety $X(M)$ and $p\colon \bar{X} \to X$ the two-fold covering map from the augmented variety to $X$. Then the following equality holds on $\bar{X}$ everywhere the right-hand term makes sense:
\begin{equation}
\label{equa:Formula1}
\tor(M,\Adr) = \frac{1}{\mT_\mu \circ p} \frac{d(1\otimes Z_\mu)}{(1 \otimes Z_\mu)}
\end{equation}
\end{proposition}
\begin{remark}
The right-hand side of Equation (\ref{equa:Formula1}) does not depend on the choice of $Z_\mu$: a direct computation shows that the right-hand side term can be written directly at any character $\chi$ where it makes sense as 
\begin{equation}
\label{equa:Formula}
\frac{2}{\mT_\mu} \frac{dI_\mu}{\sqrt{\chi(\mu)^2-4}}
\end{equation}
and neither it depends on the choice of $\mu$.
Recall that in Subsection \ref{sub:DefiTorsion} we picked the inverse convention that the one of \cite{Porti97} for the torsion, hence the term $\frac{1}{\mT_\mu}$ in (\ref{equa:Formula}). Our convention matches with the one of \cite{Dub06} but is again the inverse of \cite{DHY09}.
\end{remark}

\begin{proof}
Let $x$ be a smooth point of  $\bar{X}$, such that the function $I_\mu$ is regular at $\chi =p(x)$. Then the function $1 \otimes Z_\mu$ is regular at $x$. By Corollary \ref{cor:Smooth} the torsion can be written $f d(1 \otimes Z_\mu)$ for some invertible $f$ in $\cO_v^*$.

The differential $d(1\otimes Z_\mu)$ defines a basis element of $\det(H^1(M, \Adr))^*$ by the isomorphism $H_1(M, \Adr) \simeq \Omega^1_{k(\bar{X})/k}$. In particular its evaluation at the point $x$ defines a basis of the one dimensional vector space $\det(H^1(M, \Adbr)^*$, hence an identification $H^1(M, \Adbr) \simeq k$ given by $f \mapsto \Tr(f[\mu] \rho(\mu)_0)$ (recall that $\rho(\mu)_0$ is the trace-free matrix obtained from $\rho(\mu)$). On the other hand, in \cite{Porti97} an arbitrary generator of $H^0(\partial M, \Adbr)$ denoted by $P$ is fixed, and a basis element of $H^1(M, \Adbr)$ is given by $f \mapsto \Tr(f[\mu]P))$ (see \cite[Section 3.4]{Dub06}).

Now we picked the matrix $H$ for a basis element of $H^2(M, \Adbr)$, when in \cite{Dub06} it is again given by the matrix $P$.

In particular the choice of $P$ does not contributes to the computation of $\mT_\mu(x)$, when we have to renormalize by $\sqrt{\frac{\Tr(H^2)}{\Tr(\rho(\mu)_0^2)}}= \frac{2}{\sqrt{(\Tr(\rho(\mu)))^2-4}}$ and (\ref{equa:Formula}) follows.
\end{proof}

\begin{remark}
A consequence of Proposition \ref{prop:Formula} is that the torsion $\tor(M, \Adr)$ of this article specializes at a point $x$ of $\bar{X}$ to the usual Reidemeister torsion of the complex of $k$-vector spaces $C^*(M, \Adbr)$.
\end{remark}

\subsection{Examples}
\label{sub:Compute}
In this Subsection we compute explicit formulae for the torsion of four simple knots extriors in $\mS^3$. We will check that the degree of the divisor of the torsion (namely the sum of its zeros and poles counted with multiplicity) is the opposite of the Euler characteristic of the Riemann surface $Y$. This celebrated fact can be thought as a consequence of Riemann-Roch theorem.

\begin{example}[The trefoil knot]
Let $M$ be the exterior of the trefoil knot in $\mS^3$, with the presentation $\pi_1(M) = \langle a,b \vert \ a^2=b^3\rangle$.
Recall that the tautological representation of the component of irreducible type $X \subset X(M)$ is given by the formulae:
$$\rho(a) = \bm t&1\\-(t^2+1)&-t\ema, \rho(b)=\bm -j&0\\0 & -j^2 \ema$$
In \cite{kitano1994}, see also \cite{Dub06}, for any choice of a boundary curve $\mu \in \pi_1(\partial M)$, it is proved that $\mT_\mu(\bar{\rho})$ is a constant that does not depend on $\bar{\rho}$. Let us take for $\mu$ the meridian $ab^{-1}$, $I_\mu = (j-j^2)t$, $Z_\mu=u$, then
$$\tor(M, \Adr) =C \frac{du}{u}$$
for some non-zero constant $C$.
It has no zeros, and two poles at zero and infinity: its divisor's degree is -2. As expected, it is minus the Euler characteristic of the smooth projective model $Y$ of the augmented character variety $\bar{X}$, which is isomorphic to $\mC\mP^1$.
\end{example}

\begin{notation}
In the following examples, we will use the notation $f(t) \sim g(t)$ to say that $f$ and $g$ are equivalent around $t=0$, up to a factor that does not depends on $t$, that is $f$ and $g$ have the same vanishing order at $t=0$.
\end{notation}
\begin{example}[The figure-eight knot]
Let $M$ be the exterior of the figure-eight knot in $\mS^3$.
We take $\mu$ to be the longitude of $M$. Its trace function is $I_\mu=x^4-5x^2+2$, and $\mT_\mu(x,y)=5-2x^2$ (obtained from \cite[Example 1, p. 113]{Porti97}). From Proposition \ref{prop:Formula} it comes 
\begin{equation*}
\tor(M, \Adr) = \frac{dZ_\mu}{(5-2x^2)Z_\mu}.
\end{equation*}
A careful examination shows that it has no poles, and zeros only at infinity. We compute them now. Take $x=1/t$ a local coordinate, 
$\frac{dZ_\mu}{Z_\mu} = \frac{dI_\mu}{\sqrt{I_\mu^2-4}} \sim \frac{dt}{t}$, hence each of the four ideal points of the augmented character variety (see Subsection \ref{sub:Examples}) contribute as a zero of order 1. The divisor's degree of the torsion is $4$, and one can compute with the Riemann-Hurwitz formula that the Riemann surface $Y$ has genus $3$ (hence $\chi(Y)=-4$).
\end{example}

\begin{example}[The knot $5_2$]
\label{ex:Five}
Now $M$ is the exterior of the knot $5_2$ in $\mS^3$. Its fundamental group admits the presentation $\pi_1(M) = \langle u,v \vert \ vw=wu \rangle $ where $w=u^{-1}v^{-1}uvu^{-1}v^{-1}$.
The component of irreducible type of the character variety is described by the Riley polynomial $\phi(S,U)$ (see \cite[Section 5.2]{DHY09}). In our setting, with $x = \Tr u = \Tr v$ and $y = \Tr uv$, then $x= S^\frac{1}{2} + S^{-\frac{1}{2}}$ and $y=S+S^{-1}-U$. We obtain the following equation for the component of irreducible type of the character variety $X(M)$:
\begin{equation*}
X= \lbrace (x,y) \in k^2 \vert -x^2(y-1)(y-2)+y^3-y^2-2y+1=0 \rbrace
\end{equation*}

This affine curve admits a compactification $\hat{X}$ adding two points at infinity: an ordinary double point corresponding to the two directions $x =\infty, y=1$ or $y = 2$, and a simple point $x=y=\infty$. Apart from this, this curve is smooth. By the Noether-Pl\"ucker formula, the genus of the curve $\hat{X}$ is $g(\hat{X})=(d-1)(d-2)/2 - \de$, where $\de$ is the delta invariant. Since $d=4$ and $\de=1$, we get $g(\hat{X})=2$.

The field extension given by $\al+\al^{-1}=x$ provides a two-folds covering $Y \to \hat{X}$, that ramifies at $x^2=4$. The Riemann-Hurwitz formula implies $\chi(Y)=2\chi(\hat{X})-6=-10$, hence $Y$ is a Riemann surface of genus 6.

From \cite{DHY09} again, with $\mu$ the canonical longitude, $\mT_\mu(x,y) =  5x^4(y-2)-x^2(5y^2+7y-31)+7(y^2-y-3)$, and $I_\mu=(y^3 - 6y^2 + 12y - 8)x^{10} - (3y^4 - 10y^3 - y^2 - 68)x^8 + 3(y^5 - 43y^3 + 48y^2 + 86y - 116)x^6 + (y^6 + 6y^5 - 23y^4 - 28y^3 + 96y^2 + 28y - 105)x^4 + (2y^6 - y^5 - 16y^4 + 6y^3 + 40y^2 - 9y - 34)x^2 + 2$. 

As $\tor(M, \Adr)=\frac{dI_\mu}{\mT_\mu \sqrt{I_\mu^2-4}}$, we compute the vanishing order of the torsion at the 3 different ideal points of $\hat{X}$: 
\benu
\item
$x\sim \frac{1}{t}$, $y\sim 1+t^2$, then $\tau_\mu \sim \frac{1}{t^4}$, $\frac{dI_\mu}{\sqrt{I_\mu^2-4}}\sim \frac{dt}{t}$ and $\tor \sim t^3 dt$
\item
$x\sim \frac{1}{t}$, $y\sim 2+3t$, then $\tau_\mu \sim \frac{1}{t^2}$,  $\frac{dI_\mu}{\sqrt{I_\mu^2-4}}\sim \frac{dt}{t}$ and $\tor \sim t dt$
\item
$x \sim \frac {1}{t(1-2t^2)}$, $y \sim \frac{1}{t^2(1-2t^2)}$, then again $\tau_\mu \sim \frac{1}{t^2}$,  $\frac{dI_\mu}{\sqrt{I_\mu^2-4}}\sim \frac{dt}{t}$ and $\tor \sim tdt$
\eenu
Finally, notice that $Y \to \hat{X}$ does not ramify at infinity, hence to each ideal point of $\hat{X}$ correspond 2 ideal points of $Y$, and the divisor's degree of the differential form $\tor(M, \Adr)$ on $Y$ is 10, as expected.
\end{example}

\begin{example}[The knot $6_1$]
Here $M$ is the exterior of the knot $6_1$.
Its fundamental group admits the presentation $\pi_1(M)= \langle u,v \vert \ vw=wu \rangle $ where $w=(vu^{-1}v^{-1}u)^2$. The component of  irreducible type of the character variety is 
\begin{equation*}
X=\lbrace (x,y) \in k^2  \vert x^4(y-2)^2-x^2(y+1)(y-2)(2y-3)+(y^3-3y-1)(y-1)=0 \rbrace
\end{equation*}
The two ideal points are non ordinary double points: 
\benu
\item When $y \to 2$, $x \to \infty$, we have a double point of type $"y^2-x^6"$, its $\de$-invariant is 3.
\item When $y, x\to \infty$, we have a double point of type $"y^2-x^8"$, its $\de$-invariant is 4.
\eenu
Hence $g(\hat{X})=(d-1)(d-2)/2-\sum \de_i  =10-3-4=3$.
The covering map $Y \to \hat{X}$ given by $\al+\al^{-1}=x$ ramifies at eight finite points, thus $\chi(Y)=-16$.

When desingularizing $\hat{X}$ one obtains four ideal points, the same kind of computations as in Example \ref{ex:Five} are shortened as follows: 
\benu
\item
$x\sim \frac{1}{t(1+a t^2)}$, $y\sim \frac{2}{1+ a t^2}$ with $a$ a root of the polynomial $4z^2+6z+1$ then in both cases $\mT_\mu \sim \frac{1}{t^2}$, $\frac{dI_\mu}{I_\mu}\sim \frac{dt}{t}$ and $\tor(M, \Adr) \sim tdt$
\item
$x\sim \frac{1}{t(1-t^2)}$, $y\sim \frac{1}{t^2(1-t^2)}$, then $\mT_\mu \sim \frac{1}{t^6}$, $\frac{dI_\mu}{I_\mu}\sim \frac{dt}{t}$ and $\tor(M, \Adr) \sim t^5dt$
\item
$x\sim \frac{1}{t(1-2t^2+6t^4-25t^6)}$, $y\sim \frac{1}{t^2(1-2t^2+6t^4-25t^6)}$, then $\mT_\mu \sim 1$, $\frac{dI_\mu}{I_\mu}\sim tdt$ and $\tor(M, \Adr)\sim~tdt$
\eenu
Finally note that $Y \to \hat{X}$ is not ramified at infinity, thus the divisor's degree of $\tor(M, \Adr)$ is 16, as expected.
\end{example}
\begin{example}[The knot $7_4$, computations by S. Yoon.]
We study the manifold $M$, the exterior of the knot $7_4$. The character variety of this knot has two irreducible components that contain irreducible representations. They intersect at four finite points, that yield poles for the torsion form, as we illustrate here.

We have
$$\pi_1(M) = \langle u,v  \mid uw^2=w^2b \rangle$$ 
where $w = uv^{-1}uv^{-1}u^{-1}vu^{-1}v$. We use the coordinates fonctions $x = \Tr u = \Tr v$, and $r = \Tr uv^{-1}$, then the components of irreducible type of the character variety $X(M)$ are given by (see \cite{Chu17})
$$X_1=-1+2r^2+r^3-r^2x^2 , \qquad X_2= 1+4r-4r^2-r^3+r^4-2rx^2+3r^2x^2-r^3x^2$$
where $X_1$ is the geometric component.
The four intersection points occur when  $r= 1 \pm i$ is a root of the polynomial $r^2-2r+2$ (each giving two possible values for $x$).

In \cite{Yoon}, S. Yoon computed explicitely the torsion for twist knots, using computations from \cite{Tran}. Using this techniques, he communicated to us the following very simple formulas for the torsion form.
The torsion form on the component $X_1$ is
$$\tor_1 = \frac{2 \, \text{d}r}{(2-2r+r^2)\sqrt{1-2r^3-4r^4+r^6}}$$
and on the component $X_2$:
$$\tor_2 = \frac{2(r-1)^2 \text{d}r}{(2-2r+r^2)\sqrt{1-12r^2+42r^3-46r^4+12r^5+9r^6-6r^7+r^8}}$$
One can check that the two torsions forms have four poles (of order one) at the intersections points $(r^2-2r+2)=0$, and nowhere else.
\label{ex:74}
\end{example}
\section{Ideal points, torsion form and essential surfaces}
\label{sec:Ideal}
In this section we study the behavior of the torsion form at ideal points of the augmented character variety. Ideal points have been shown in Subsection \ref{sub:CS} to produce essential surfaces in the manifold $M$. In what follows we will do the following assumption on the essential surfaces produced as such. We fix $v$ an ideal point in the smooth projective model $Y$, the essential surface $\Sigma$ dual to $v$ will be assumed to be:
\benu
\item
union of $n$ parallel copies $\Sigma_i$
\item
separating, in the sense that for any $i$, $M \setminus \Sigma_i = M_1\cup M_2$ is not connected
\item
free, that is the connected components $M_1$ and $M_2$ of $M \setminus \Sigma_i$ are handlebodies
\eenu
By Lemma \ref{lem:Converge} the divergent tautological representation $\rho \colon \pi_1(M) \to \SL_2(k(Y))$ restricts to a convergent representation $\rho_\Sigma\colon \pi_1(\Sigma) \to \SL_2(\cO_v)$, such that the residual representation $\brho_\Sigma \colon \pi_1(\Sigma) \to \SL_2(k)$ is reducible.

In what follows we will assume that the representation $\rho_\Sigma$ is irreducible and not residually central.
The aim of this section is to prove the following theorem:
\begin{theorem}
\label{theo:Ideal}
Let $\bar{X}$ be a one dimensional component of irreducible type, essentially reduced in the augmented variety $\bar{X}(M)$. Let $v$ be an ideal point of the smooth projective model $Y$ of $\bar{X}$, and $\Sigma$ a dual free separating essential surface in $M$, that is union of $n$ parallel copies $\Sigma_i$. Assume furthermore that the restricted representation $\rho_\Sigma$ is irreducible and not residually central. Then the following inequality holds:
\begin{equation*}
v(\tor(M, \Adr)) \le -n(\chi(\Sigma_i)+1)
\end{equation*}
\end{theorem}

Combining this theorem with the results of Section \ref{sec:Finite}, we obtain the following corollary:
\begin{corollary}
Let $M$ be a 3-manifold with rational homology of a circle, whose Alexander polynomial $\Delta_M$ has only simple roots. Let $\bar{X}$ be a smooth one dimensional component of irreducible type, essentially reduced in the augmented variety $\bar{X}(M)$, such that to each ideal point $v$ of the smooth projective model $Y$ of $\bar{X}$, it can be associated a dual surface $\Sigma_v$ that satisfies the hypotheses of Theorem \ref{theo:Ideal}. Then we have
\begin{equation*}
\chi(Y)\ge \sum_{v \text{ ideal }} n_v(\chi(\Sigma_v)+1)
\end{equation*}
\end{corollary}
\begin{proof}
Denote by $\rho\colon \pi_1(M) \to \SL_2(k(Y))$ the tautological representation.
Since $\Delta_M$ has simple roots and $\bar{X}$ is smooth, the torsion $\tor(M, \Adr)$ does not vanish at finite points (see Theorem \ref{theo:Singular}). Hence by Theorem \ref{theo:Ideal}, its divisor's degree is bounded by $\sum\limits_{v \text{ ideal }} -n_v(\chi(\Sigma_v)+1)$, now it follows from Riemann-Roch Theorem the divisor's degree of a differential form on a Riemann surface $Y$ is equal to minus the Euler characteristic of $Y$.
\end{proof}

\begin{remark}
\benu
\item
The hypothesis that the dual surface $\Sigma$ is separating excludes Seifert surfaces, but it is known that a separating essential surface can always be produced by the Culler-Shalen theory at some ideal point (it is the way the weak Neuwirth conjecture is proven in \cite{CS84}). It also known that if the Seifert surface is a fiber, then it cannot be the dual surface of an ideal point.
\item 
If $M$ is irreducible and small (does not contain any closed essential surfaces), then any essential surface is free (see \cite[Proposition 4.3]{Pri}). In general the assumption that $\Sigma$ is free exclude closed essential surfaces.
\item
If $\bar{X}$ covers the character of a faithful representation, then $\rho_\Sigma$ is irreducible: if not the commutator subgroup of $\pi_1(\Sigma)$ is faithfully mapped onto an abelian group, a contradiction. 
\eenu
\end{remark}

The section is organized as follows:  in Subsection \ref{sub:IdealExamples} we compare the inequality of Theorem \ref{theo:Ideal} with the examples of Section \ref{sub:Compute}, then in Subsection \ref{sub:ProofIdeal} we prove Theorem \ref{theo:Ideal}.

\subsection{Examples}
\label{sub:IdealExamples}

Incompressible surfaces of two-bridge knots are classified in \cite{HT85}.
\subsubsection{The trefoil knot.}

The essential dual surface $\Sigma$ is an annulus, hence $\rho_\Sigma: \mZ \to \SL_2(k(t))$ is abelian, and the theorem cannot apply (nevertheless the torsion has a pole of order $1=-\chi(\Sigma)-1$ at the ideal points corresponding to $\Sigma$).
 
 \subsubsection{The figure-eight knot.}

There are two essential dual surfaces $\Sigma_1$ and $\Sigma_2$ that are two-holed tori, and the torsion vanishes at order $1$ at each ideal point. The inequality of Theorem \ref{theo:Ideal} is an equality $1 = -\chi(\Sigma_i)-1$, with $n=1$. 
\subsubsection{ The knot $5_2$.}

There are two separating essential surfaces $\Sigma_1$ and $\Sigma_2$, see Figure \ref{noeud}, and as explained in the introduction, again the equality of Theorem \ref{theo:Ideal} holds.
\subsubsection{The knot $6_1$.}

Again, there are two separating essential surfaces, the first one has Euler characteristic equal to $-2$ (a two-holed torus), and the second one has Euler characteristic equal to $-6$ (a two-holed genus 3 surface). The third essential surface is the Seifert surface. Those three surfaces are detected by the character variety. At those ideal points the vanishing order of $\tor(M)$ is $1$ (three times) and $5$, the latter corresponds to the genus 3 non-Seifert surface.
\begin{remark}
Those examples give an insight of the power of Theorem \ref{theo:Ideal}: essential surfaces of those knots are completely understood, but we can even deduce informations on the Bass-Serre tree constructed at ideal points of the character variety.
\end{remark}

\subsection{Proof of Theorem \ref{theo:Ideal}}
\label{sub:ProofIdeal}

In this subsection we use the notations of Lemma \ref{lem:Converge}. We fix a connected component $\Sigma_1$ of $\Sigma$, that we will denote by $\Sigma$ to avoid heavy notations, and we can forget about the other components and think that the dual surface is connected, excepted that the matrix $U_n=\bsm t^n & 0\\0 & 1\esm$ of Lemma \ref{lem:Converge} carries the number $n$ of connected components. Recall that $M = M_1 \cup_{\Sigma} M_2$, with $M_1$ and $M_2$ handlebodies. The representation $\rho_1\colon \pi_1(M_1) \to \SL_2(\cO_v)$ is convergent, $\rho_2\colon \pi_1(M_2) \to \SL_2(k(Y))$ is not but $\rho'_2 = U_n^{-1} \rho_2 U_n$ is.
 
There is an exact sequence of complexes of $k(Y)$-vector spaces
\begin{equation}
\label{equa:Mayer}
0 \to C^*(M, \Adr) \to C^*(M_1, {\Adr}_1)\oplus C^*(M_2, {\Adr}_2) \to C^*(\Sigma, {\Adr}_\Sigma) \to 0
\end{equation}
Since $\rho_\Sigma$ is not abelian, the vector spaces $H^0(\Sigma, {\Adr}_\Sigma)$ and $H^0(M_i, {\Adr}_i)$ vanish
and the splitting (\ref{equa:Mayer}) induces the long exact sequence in cohomology
\begin{align}
\label{equa:ExactSeq1}
0\to H^1(M, {\Adr}) \to H^1(M_1, {\Adr}_1) \oplus H^1(M_2, {\Adr}_2)\\ \nonumber \to H^1(\Sigma, {\Adr}_\Sigma) \to H^2(M, \Adr) \to 0
\end{align}
that we denote by $\mathcal H_1$.

We will use the multiplicativity formula (see Proposition \ref{prop:Multiplicativity}):
\begin{equation}
\label{equa:Milnor}
\tor(M, \Adr) =\frac{ \tor(M_1, {\Adr}_1)\tor(M_2, {\Adr_2})}{\tor(\Sigma, {\Adr}_\Sigma)} \tor(\mathcal H_1)
\end{equation}
In a first step, we need to make sense of all the factors involved in (\ref{equa:Milnor}): we pick geometric bases $\textbf{c}_1,\textbf{c}_2$ and $\textbf{c}_\Sigma$ of the complexes $C^*(M_1, {\Adr}_1)$, $C^*(M_2, {\Adr}_2)$, and $C^*(\Sigma, {\Adr}_\Sigma)$, and arbitrary homological basis $\textbf{h}_1$, $\textbf{h}_2$ and $\textbf{h}_\Sigma$ of their homology.
Hence the torsions of those complex are defined as in \ref{subsub:Definition} as elements of $k(Y)^*$. Now the sequence (\ref{equa:ExactSeq1}) is thought as an acyclic complex, and its torsion depends on the choices $\textbf{h}_1$, $\textbf{h}_2$ and $\textbf{h}_\Sigma$, and of a basis element in $\det(H^*(M, \Adr))$. As we did in Section \ref{sub:ProofFinite}, for the latter we fix the element $dt \otimes H^*$.

Since $\rho_1$, $\rho'_2$ and $\rho_\Sigma$ are convergent, we can define the complexes of $\cO_v$-modules $C^*(M_1, {\Adr}_1)_v$, $C^*(M_2, {\Adr}_2)_v$, $C^*(\Sigma, {\Adr}_\Sigma)_v$. We compute their cohomology as $\cO_v$-modules, and their residual cohomology:
\begin{lemma}
\label{lem:ResidualSigma}
For $i=0$ or $i\ge 2$, the $\cO_v$-modules $H^i(M_1, {\Adr}_1)_v$, $H^i(M_2, {\Adr}'_2)_v$ and $H^i(\Sigma, {\Adr}_\Sigma)_v$ are trivial, as well as the corresponding residual $k$-vector spaces.

Furthermore $H^1(M_1, {\Adr}_1)$ and $H^2(M_2, {\Adr}'_2)$ are free of rank $-\frac 3 2 \chi(\Sigma)$ and $H^1(\Sigma, {\Adr}_\Sigma)$ is free of rank $-3 \chi(\Sigma)$. The corresponding residual $k$-vector spaces have the same dimensions.
\end{lemma}
\begin{proof}
First notice that those manifold have the same homotopy type as a graph, hence the statement for $i\ge2$ is clear.

The representation $\rho_\Sigma$ is irreducible, hence $\rho_1$ and $\rho'_2$ also. It implies that $H^0(\Sigma, {\Adr}_\Sigma)_v$, $H^0(M_1, {\Adr}_1)_v$ and $H^0(M_2, {\Adr}_2)_v$ are trivial.

By hypothesis $\rho_\Sigma$ is not residually central, hence by Lemma \ref{lem:NotAbelian} one can chose the tautological representation such that $\rho_\Sigma$ is not residually abelian. Thus we have $H^0(\Sigma, {\Adbr}_\Sigma)= \lbrace 0 \rbrace$. Since $\brho_\Sigma$ is not abelian, neither is $\brho_1$, and a careful examination of the proof of Lemma \ref{lem:NotAbelian} shows that $\brho_2$ neither, hence $H^0(M_1, {\Adbr}_1)$ and $H^0(M_2, {\Adbr}_2)$ are trivial.

An argument such as in Subsection \ref{sub:Twist} shows now that $H^1(\Sigma, {\Adr}_\Sigma)_v, \ldots$ are free $\cO_v$-modules, with rank prescribed by the Euler characteristic of $\Sigma$, and the last statement is then clear.
\end{proof}
\begin{proposition}
\label{prop:TorsionInvertible}
There is a choice of bases $\textbf{h}_1$, $\textbf{h}_2$ and $\textbf{h}_\Sigma$ such that the factors $\tor(M_1,{\Adr}_1)$, $\tor(M_2, {\Adr_2})$ and $\tor(\Sigma, {\Adr}_\Sigma)$ of (\ref{equa:Milnor}) are elements of $\cO_v^*$. In particular the vanishing order of $\tor(M, \Adr)$ can be computed as the vanishing order of $\tor(\mathcal H_1, \textbf{h}_1, \textbf{h}_2, \textbf{h}_\Sigma)$, the torsion of the Mayer--Vietoris long exact sequence (\ref{equa:ExactSeq1}).
\end{proposition}

\begin{proof}
Those factors are torsion of based complex of $k(Y)$-vector spaces with based homology, hence they lie in $k(Y)^*$ by definition. 
In addition, we can chose the bases $\textbf{c}_\Sigma, \textbf{c}_1,\textbf{c}_2,\textbf{h}_\Sigma, \textbf{h}_1, \textbf{h}_2$ of the paragraph above such that they provide a generating set for the corresponding terms as $\cO_v$-modules, because the $\cO_v$-modules $C^*(\Sigma, \Adr_\Sigma)_v, \ldots, H^*(M_2, \Adr'_2)_v$ are free, and those choices do not affect the computation of the torsion. 
To be precise, we assume that we have chosen for instance a basis $h_2$ of the free $\cO_v$-module $H^1(M_2, \Adr'_2)_v$ that spans $H^1(M_2, \Adr'_2)$ as a $k(Y)$-vector space, and that it is mapped on a basis through the isomorphism of $k(Y)$-vector spaces $H^1(M_2, \Adr'_2) \to H^1(M_2, \Adr_2)$.  Finally, the map $H^1(M_1, \Adr_1)_v \to H^1(\Sigma, \Adr_\Sigma)_v$ identifies the basis $h_1$ to a sub-basis of $h_\Sigma$.

Now we prove that the torsions of those complexes lie in $\cO_v^*$: let us perform the computation for, say, $M_1$.
The complex is $C^0(M_1, \Adr_1)_v \xrightarrow{A} C^1(M_1, \Adr_1)_v$. Since $H^0(M_1, \Adr_1)_v$ is trivial, the matrix $A$ is the matrix of an injective $\cO_v$-linear morphism. Moreover, $H^1(M_1, \Adr_1)_v$ is free, hence the determinant of the restriction $\bar{A}: C^0(M_1, \Adr_1)_v \to \operatorname{im}(A)$ is invertible as claimed.

Since $\tor(M_1), \tor(M_2)$ and $\tor(\Sigma)$ take values in $\cO_v^*$,  the valuation of $\tor(M, \Adr)$ is determined by the torsion of the exact sequence $\mathcal H_1$.
\end{proof}

Now in order to compute the torsion of the sequence $\mathcal H_1$, we modify it slightly.
The representation $\rho_2$ is conjugated to a convergent representation $\rho_2' \colon \pi_1(M_2) \to \SL_2(\cO_v)$ by the matrix $U_n=~\bsm t^n & 0\\ 0 & 1 \esm$. It yields an isomorphism $H^1(M_2, \Ad\circ{\rho'_2}) \xrightarrow{\sim} H^1(M_2, \Ad\circ{\rho_2})$ given by $\zeta_2 \mapsto U_n \zeta_2 U_n^{-1}$. Hence we can rewrite the sequence $\mathcal H_1$ as a new exact sequence $\mathcal H_2$ of $k(Y)$-vector spaces, given by:
\begin{align}
\label{equa:ExactSeq2}
0 \to H^1(M, \Adr) \xrightarrow{d_1} H^1(M_1, \Ad\circ{\rho_1}) \oplus H^1(M_2,\Ad\circ{\rho'_2}) \\ \nonumber \xrightarrow{d_2} H^1(\Sigma, \Ad\circ{\rho_\Sigma})\xrightarrow{\de} H^2(M, \Adr) \to 0
\end{align}

Under the isomorphism $\alpha_* \colon H^1(M, \Adr) \simeq (\Omega^1_{k(Y)/k})^*$ described in (\ref{equa:LongSequence}), the basis element $dt \in \det(H^1(M, \Adr))^*$ corresponds to a derivation that we denote by $\frac d{dt} \colon k(Y) \to k(Y)$ (see Proof of Lemma \ref{lem:Onto1}). We keep the same notation for the induced map 
$$\frac d{dt} \colon \SL_2(k(Y))\to \SL_2(k(Y).$$
It yields an explicit basis element $(\frac d{dt}\rho) \rho^{-1} \in H^1(M, \Adr)$.
In the following lemma we compute the morphisms $d_1, d_2$ and $\delta$ involved in the sequence $\mathcal H_2$ introduced in (\ref{equa:ExactSeq2}).
\begin{lemma}
\label{lem:Morph}
The morphisms $d_1, d_2$ and $\delta$ are
$$d_1\colon (\frac{d}{dt}\rho)\rho^{-1} \mapsto \left((\frac{d}{dt}\rho_1)\rho_1^{-1}, (\frac{d}{dt}\rho'_2){\rho'_2}^{-1}\right),$$
$$d_2\colon(\zeta_1, \zeta_2)  \mapsto  \zeta_1- U_n \zeta_2 U_n^{-1}$$
and  
$$\de\colon \xi \mapsto \Tr\left(\xi([\partial \Sigma]) H\right)$$
with the natural identification $H^2(M, \Adr) \simeq k(Y)$ described in (\ref{equa:NaturalIsom}).
\end{lemma}

\begin{proof}
The morphism $d_1$ is given by:
$$d_1((\frac{d}{dt}\rho)\rho^{-1}) =  \left((\frac{d}{dt}\rho_1)\rho_1^{-1},  U_n^{-1}  (\frac{d}{dt}\rho_2)\rho_2^{-1} U_n\right)$$
and we have
\begin{align*}
 U_n^{-1}  (\frac{d}{dt}\rho_2)\rho_2^{-1} U_n&=U_n^{-1}\frac{d}{dt}(U_n\rho'_2U_n^{-1}) (U_n{\rho'_2}^{-1}U_n^{-1}) U_n\\
 &=(\frac{d}{dt}\rho'_2 + U_n^{-1}\frac{d}{dt}U_n \rho_2'+\rho_2' \frac{dU_n^{-1}}{dt} U_n){\rho'_2}^{-1}\\
 &=(\frac{d}{dt}\rho'_2){\rho'_2}^{-1} +\frac{1}{t^n} (\rho'_2 N {\rho'_2}^{-1}-N)\\
 \end{align*}
 where the matrix $N$ is $\bsm 1&0\\0&0 \esm$.
 Now the term $(\rho'_2 N {\rho'_2}^{-1}-N)$ is the coboundary of the element $N_0=N-\frac{1}{2}I$ in $C^1(M_2, \Adr'_2)$, hence it vanishes in $H^1(M_2, \Adr'_2)$ and the claimed formula for $d_1$ follows.
 The formula for $d_2$ is just the composition of the usual morphism in Mayer--Vietoris sequences, composed with the isomorphism $H^1(M_2, \Ad\circ{\rho'_2}) \xrightarrow{\sim} H^1(M_2, \Ad\circ{\rho_2})$.
 
For the third morphism $\delta$, the naturality of the Mayer-Vietoris sequence and the exact sequence of a pair yield the commutative diagram.
\begin{center}
\begin{tikzpicture}
\tikzset{node distance=2.5cm, auto}
\node(A) {$H^1(\Sigma, \Ad\circ{\rho})$};
\node(B)[right of=A, xshift=1cm] {$H^2(M, \Ad\circ\rho)$};
\node(C)[below of=A, yshift=1cm] {$H^1(\partial \Sigma, \Ad\circ\rho)$};
\node(D)[right of=C, xshift=1cm] {$H^2(\partial M, \Ad\circ\rho)$};
\node(E)[right of=D, xshift=1cm] {$H^0(\partial M, \Ad\circ\rho)^*$};
\draw[->] (A) to node{$\de$} (B);
\draw[->] (A) to (C);
\draw[->] (C) to (D);
\draw[->] (B) to (D);
\draw[->] (D) to node{$\sim$}(E);
\end{tikzpicture}
\end{center}
As the second vertical arrow is an isomorphism, it is enough to compute the composition of the maps $H^1(\Sigma, \Ad\circ{\rho}) \to H^1(\partial \Sigma, \Ad\circ\rho) \to H^2(\partial M, \Ad\circ\rho) \to k(Y)$, which is simply $\xi \mapsto \Tr(\xi([\partial \Sigma])H)$, see (\ref{equa:NaturalIsom}).
\end{proof}

Now each term of the sequence $\mathcal H_2$ given in (\ref{equa:ExactSeq2}) can be thought as an $\cO_v$-module tensorized by $k(Y)$, but the map $d_2$ above does not restrict to a morphism of $\cO_v$-modules. Hence we will consider the following exact sequence, that we denote by $\mathcal H_3$:
\begin{equation}
\label{equa:ExactSeq3}
 0 \to k(Y) \xrightarrow{d_1} H^1(M_1, \Ad\circ{\rho_1}) \oplus H^1(M_2,\Ad\circ{\rho'_2}) \xrightarrow{t^n \cdot d_2} H^1(\Sigma, \Ad\circ{\rho_\Sigma}) \xrightarrow{\de} k(Y) \to 0
 \end{equation}
where we just have multiplied the morphism $d_2$ by $t^n$. We will denote by $D_2$ this new map, which restricts to  morphism of $\cO_v$-modules $H^1(M_1, \Ad\circ{\rho_1})_v \oplus H^1(M_2,\Ad\circ{\rho'_2})_v \xrightarrow{D_2} H^1(\Sigma, \Ad\circ{\rho_\Sigma})_v$.
Hence the sequence $\mathcal H_3$ can be seen as the following sequence $\mathcal H_4$ tensored by $k(Y)$:
\begin{equation}
\label{equa:H4}
 0 \to \cO_v \xrightarrow{d_1} H^1(M_1, \Ad\circ{\rho_1})_v \oplus H^1(M_2,\Ad\circ{\rho'_2})_v \xrightarrow{D_2} H^1(\Sigma, \Ad\circ{\rho_\Sigma}v \xrightarrow{\de} \cO_v \to 0
\end{equation}

From now on we suppose that the choices of bases we made $h_\Sigma, h_1$ and $h_2$ gave splittings $H^1(M_1, \Ad\circ{\rho_1}) \oplus H^1(M_2,\Ad\circ{\rho'_2})= \ker d_2 \oplus E_1$, and $H^1(\Sigma, \Ad\circ{\rho_\Sigma}) = d_2(E_1) \oplus E_2$. Let $\Delta_2$ be the restricted map ${D_2}_{\vert_{E_1}}: E_1 \to d_2(E_1)$. 
\begin{lemma} 
\label{lem:Same}
The torsion of the exact sequence $\mathcal H_3$ given in (\ref{equa:ExactSeq3}) is
\begin{equation*}
\tor(\mathcal H_3) = \frac{1}{\det \Delta_2} c \text{, with } c \in \cO_v^*
\end{equation*}
\end{lemma}

\begin{proof}
Consider the definition of the torsion of \ref{subsub:Cayley}, the following equality holds:
\begin{equation*}
\tor(\mathcal H_3)=\frac{ \det(d_1: k(Y) \to d_1(k(Y))) \det(\de: E_2 \to k(Y)) }{\det \Delta_2}
\end{equation*}
Then we conclude the proof by noting that the numerator lies in $\cO_v^*$, since by Lemma \ref{lem:Morph} neither the morphism $d_1$ nor $\delta$ involve the parameter $t$.
\end{proof}

Hence we are now reduced to compute $v(\det(\Delta_2))$. 
To do this, the idea is the following: recall that the completion of the valuation ring $\cO_v$ is isomorphic to $k[[t]]$, the ring of formal series.
Consider a matrix $A \in \mathcal{M}_n (\cO_v)$ as a formal series $A=\sum t^i A_i$, with $A_i \in \mathcal{M}_n(k)$, the question is to compute the valuation of its determinant (namely, its vanishing order when $t=0$). 

If $\det A_0 \neq 0$, then $A$ is invertible, $\det A \in \cO_v^*$ and $v(\det A)=0$. If not, we have $k^n \xrightarrow{A_0}  k^n$ which is not invertible and define $H^0(A_0)=\ker A_0$, $H^1(A_0) = \coker A_0$. Both are non-trivial $k$-vector space of the same dimension.
Pick $P,Q \in \operatorname{GL}_n(k)$ such that $P A_0 Q = \bsm 0 & 0\\0 & \ I_{n-r_0} \esm $ is diagonal, where  $r_0= \dim \ker A_0$, and $I_{n-r_0}$ is the $(n-r_0)$ identity matrix. Then to compute $\det A$, we need to compute the determinant of the $r_0 \times r_0$ first block of $A_1 + tA_2 ...$. More precisely $\det A=t^{r_0} \det A'_1 + o(t^{r_0})$, where $A'_1$ is the restriction of $\sum t^i A_{i+1}$ to $H^0(A_0)\otimes k[[t]]$, followed by the projection $k[[t]]^n \to H^1(A_0) \otimes k[[t]]$.

One proceeds by induction, the argument is formalized in Lemma \ref{lem:Determinant}, but before that we fix some notations.
Let $A$ be a matrix in $\mathcal{M}_n (\cO_v)$ such that $\det A \neq 0$. Working in the completion $\hat{\cO}_v$ if necessary, we define $A_{\geq 0}=A$, and inductively $A_{\geq i +1} = \frac{d}{dt} A_{\geq i}$ restricted to $H^0( A_{\geq i}(0) )\otimes k[[t]]$ followed by the projection   $k[[t]]^{n-\sum_{k=0}^i r_k} \to H^1(A_{\geq i}(0)) \otimes k[[t]]$, where $r_i = \dim \ker A_{\geq i}(0)$.

\begin{lemma} 
\label{lem:Determinant}
There is a constant $c\in \cO_v^*$ such that $\det(A)=t^{\sum r_i} c$, in particular $v(\det(A)) = \sum r_i$.
\end{lemma}

\begin{proof}
 As $\det A \neq 0$, there is an $i_0$ such that $r_{i_0}= 0$ and $r_i \neq 0$ for $i<i_0$. If $i_0=0$ then $\det A \in \cO_v^*$ and we are done. Else, take $0 < i \leq i_0$. After fixing appropriated bases of $\ker A_{\geq i-1}(0)$, one can write $A_{\geq i}(0)$ as a diagonal matrix $\bsm 0 & 0\\0 & \ I_{n-r_{i-1}} \esm$.
 
Developing the determinant of $A_{\geq i}$ along its rows, one gets that $\det A_{\geq i} = t^{r_i} \det A_{\geq i+1} + o(t^{r_i})$, and the result follows by induction.
\end{proof}

We will apply Lemma \ref{lem:Determinant} to the morphism $\Delta_2$ of Lemma \ref{lem:Same}.
Recall that for each convergent representations $\rho_1, \rho_2',\rho_\Sigma$, we have the residual representations $\bar{\rho}_1, \bar{\rho}_2, \bar{\rho}_\Sigma$ taking values in the residual field $k$. Moreover, $\brho_\Sigma = \brho_{1, \Sigma}=\overline{U_n \rho'_{2, \Sigma} U_n^{-1}}$ is reducible, non abelian, thus we have:
\begin{lemma}
The residual representations have the form
\begin{equation*}
\brho_{1,\Sigma} = \bm \la & 0\\ \la u_1 &\la^{-1} \ema, \brho_{2,\Sigma} = \bm \la & \la^{-1}u_2 \\ 0 & \la^{-1} \ema
\end{equation*}
with $\la \in H^1(\Sigma,k^*)$ a homomorphism and  $u_1 \in H^1(\Sigma,\la^{-2}), u_2 \in H^1(\Sigma,{\la^{2}})$ non trivial cocycles.
\end{lemma}
\begin{proof}
The expression of $\brho_{2,\Sigma}$ follows from the conjugacy formula $\rho_1 = U_n\rho'_2U_n^{-1}$ when restricted on $\pi_1(\Sigma)$, the $u_i$'s are non trivial because the residuals representations are not abelian, see Lemma~\ref{lem:Cocycleu}.
\end{proof}
From (\ref{equa:H4}) we obtain the (non-acyclic) complex of $k$-vector spaces $\mathcal H_5$:
\begin{equation}
\label{equa:ExactSeq4}
0 \to k \xrightarrow{\bar{d}_1} H^1(M_1, \Ad\circ{\brho_1}) \oplus H^1(M_2, \Ad\circ{\brho_2}) \xrightarrow{\bar{D}_2} H^1(\Sigma, \Ad\circ{\brho_{1,\Sigma}}) \xrightarrow{\bar{\de}} k \to 0
\end{equation}
with $\bar{d}_1(1)= (v_1,v_2)$ and $\bar{D}_2(\zeta_1, \zeta_2) = z_{2,\Si}F$, where $z_{2,\Sigma}$ denotes the lower-left entry of $\zeta_2$, restricted to $\pi_1(\Sigma)$.

As in Subsection \ref{subsub:step1}), the triangularity of the adjoint action of $\brho_{i,\Sigma}$ provides the following splittings:
\begin{align*}
0 \to k^2_{(\Adbr_{1,\Sigma})_{3,3}} \to &\slf_2(k)_{\Ad\circ{\brho_{1,\Sigma}}} \to \varphi_{\la^{2}} \to 0\\
0 \to k^2_{(\Adbr_{2,\Sigma})_{3,3}} \to &\slf_2(k)_{\Ad\circ{\brho_{2,\Sigma}}} \to \varphi_{\lambda^{-2}} \to 0\\
0 \to \varphi_{\lambda^{-2}} \to &k^2_{(\Adbr_{1,\Sigma})_{3,3}} \to k \to 0
\end{align*}
and thus the exact sequences of $k$-vector spaces:
\begin{align}
\label{37}
0 \to H^1(\Sigma, (\Adbr_{2,\Sigma})_{3,3}) \to &H^1(\Sigma, \Ad\circ{\brho_{2,\Sigma}}) \xrightarrow{p} H^1(\Sigma, \varphi_{\la^{-2}}) \to 0\\ 
\label{38}
0 \to H^0(\Sigma) \to H^1(\Sigma, \varphi_{\la^{-2}}) \to &H^1(\Sigma, (\Adbr_{1,\Sigma})_{3,3}) \to H^1(\Sigma) \to 0\\
\label{39}
0 \to H^1(\Sigma, (\Adbr_{1,\Sigma})_{3,3}) \to &H^1(\Sigma, \Ad\circ{\brho_{1,\Sigma}}) \to \ldots
\end{align}
We denote by $j$ the composition 
\begin{align}
\label{j}
j \colon H^1(\Sigma, \varphi_{\la^{-2}}) \to H^1(\Sigma, (\Adbr_{1,\Sigma})_{3,3})\to H^1(\Sigma, \Ad\circ{\brho_{1,\Sigma}})
\end{align}
 of (\ref{38}) and (\ref{39}).
\begin{lemma}
\label{lem:Kerj}
The space $\ker j$ is one dimensional. It is generated by the image of $H^0(\Sigma)$ in $H^1(\Sigma,\la^{-2})$, that is by $\partial_{1,\Sigma}  H =-2u_1F$.
\end{lemma}
\begin{proof}
The first statement is clear since $H^0(\Sigma)$ is one-dimensional. All we need to do is to compute $\partial_{1,\Sigma} H = \brho_{1,\Sigma} H {\brho_{1,\Sigma}}^{-1} -H$, and we obtain the claimed result. 
\end{proof}

The inclusion $\Sigma\subset M_2$ induces $i: H^1(M_2, \Ad\circ{\brho_2}) \to H^1(\Sigma, \Ad\circ{\brho_{2,\Sigma}})$.
Let $K_{j \circ p \circ i}$ be the kernel of the morphism $j \circ p \circ i: H^1(M_2, \Ad\circ{\brho_2}) \to H^1(\Sigma, \Ad\circ{\brho_{1,\Sigma}})$ obtained by composing the morphism $i$ above with the morphism $p$ of (\ref{37}) and the morphism $j$ of (\ref{j}).
\begin{lemma}
\label{lem:DimKer}
We have
 $$\dim(K_{j \circ p \circ i}) \ge -\frac{1}2 \chi(\Sigma)+1.$$
\end{lemma}
\begin{proof}
We know that 
\begin{align*}
\dim H^1(\Sigma, \Ad\circ{\brho_\Sigma})=-3\chi(\Sigma),\\  \dim H^1(\Sigma, \varphi_{\la^2})= -\chi(\Sigma), \\ \dim H^1(M_i, \Ad\circ{\brho_i})=~-\frac3 2 \chi(\Sigma).
\end{align*}
Since the morphism $p$ is onto, the dimension of the space $\ker p $ is 
\begin{equation*}
\dim(\ker p) = \dim H^1(\Sigma, \Ad\circ{\brho_{2,\Sigma}}) - \dim H^1(\Sigma, \varphi_{\la^2}) = -3\chi(\Sigma) - (-\chi(\Sigma))=-2\chi(\Sigma).
\end{equation*}
 If $i$ is injective and if $\ker p $ and $\im i $ intersect transversally, then $\dim \ker p \cap \im i =  -\frac 1 2\chi(\Sigma)$. We define the integer $s$ by the formula $\dim \ker p \cap \im i = -\frac 1 2 \chi(\Sigma) +s$. By Lemma \ref{lem:Kerj} the space $\ker j$ has dimension 1, moreover, $p \circ i (v_2) =u_1$ is non-trivial in $\ker j$, hence the dimension of the intersection $\left(\ker (j \circ p) \cap  \im i\right)$ is equal to $-\frac{\chi(\Sigma)}2+1+s$. Now if $i$ is injective, the latter is the dimension of the space $\ker(j\circ p\circ i)$ while if $i$ has non trivial kernel, it will possibly increase this dimension. In any case the inequality 
 \begin{equation*}
 \dim (\ker (j\circ p\circ i)) \ge -\frac1 2 \chi(\Sigma) +1
 \end{equation*}
 holds.
\end{proof}
Now we can prove Theorem \ref{theo:Ideal}:
\begin{proof}
We use the notations of Lemma \ref{lem:Determinant}. We compute $r_0$, the dimension of the first homology group of the sequence $\mathcal H_5$ of (\ref{equa:ExactSeq4}), namely $H^1(\mathcal H_5)=\ker \bar{D}_2/ \im \bar{d}_1$. By Lemma \ref{lem:DimKer} we have $$\dim \ker \bar{D}_2 \ge -\chi(\Sigma)/2+1 + (-3/2 \chi(\Sigma))=-2 \chi(\Sigma) +1.$$
Hence $r_0 \ge -2 \chi(\Sigma) $. 

We use notations of Lemma \ref{lem:Determinant}: observe that the higher order maps $\frac{d^i}{d^it} \vert_{t=0}(D_2)$, for $i = 1, \ldots, n-1$, remain zero when restricted to $\ker(\bar{D}_2)$ (because of the factor $t^n$ in the map $D_2$), hence each $r_i, i=1, \ldots, n-1$ is greater than $-2 \chi(\Si)$.
Let $r = \sum_{i\geq n} r_i$, we have from Lemma \ref{lem:Determinant} that $$\det(D_2)= \sum_{i\ge 0} r_i \ge -2 n \chi(\Sigma) +r,$$ hence $v(\tor(\mathcal H_3)) \le 2n \chi(\Sigma)$ by Lemma \ref{lem:Same}.

Finally, it is easy to obtain the following relation between the torsions of the sequence $\mathcal H_2$ of (\ref{equa:ExactSeq2}) and of the sequence $\mathcal H_3$ of (\ref{equa:ExactSeq3}):
\begin{equation*}
t^{n \rk(d_2)} \tor(\mathcal H_3) = \tor(\mathcal H_2) 
\end{equation*}. Since $\rk(d_2) = -3 \chi(\Sigma) -1$,  we obtain 
\begin{equation*}
v(\tor(\mathcal H_2) \le 2n\chi(\Sigma) -3n \chi(\Sigma) -n = -n(\chi(\Sigma)+1)
\end{equation*}
and the theorem follows now from Proposition \ref{prop:TorsionInvertible}.
\end{proof}

%
%

\bibliographystyle{plain} 
\bibliography{biblio} 
\end{document}